\documentclass[12pt,a4paper, reqno]{amsart}
\usepackage[utf8]{inputenc}
\usepackage{mathptmx}
\usepackage{amssymb}
\usepackage{hyperref}
\usepackage{mathrsfs}
\usepackage[mathscr]{euscript}
\usepackage{booktabs}

\usepackage{todonotes}

\usepackage[ignoreunlbld,norefs,nocites,nomsgs]{refcheck}

\usepackage{pgf, tikz}

\usepackage{fancyvrb} \RecustomVerbatimEnvironment{Verbatim}{Verbatim}
{xleftmargin=15pt, frame=single, fontsize=\small}

\usepackage{color}

\definecolor{light}{gray}{0.9}
\definecolor{medium}{gray}{0.8}

\newtheorem{theorem}{Theorem}[section]
\newtheorem{lemma}[theorem]{Lemma}
\newtheorem{corollary}[theorem]{Corollary}
\newtheorem{proposition}[theorem]{Proposition}

\theoremstyle{definition}
\newtheorem{remark}[theorem]{Remark}
\newtheorem{definition}[theorem]{Definition}

\numberwithin{equation}{section}

\def\ZZ{{\mathbb Z}}
\def\NN{{\mathbb N}}

\def\11{{\mathbb 1}}
\def\cR{{\mathcal R}}

\def\cU{{\mathcal U}}

\def\cB{{\mathcal B}}

\def\cF{{\mathcal F}}
\def\cM{{\mathcal M}}

\def\FV{FV}
\def\New{{\mathcal N}}

\def\chara{\operatorname{char}}

\def\iso{\cong}

\def\Rees{\operatorname{\cR}}

\let\epsilon=\varepsilon

\allowdisplaybreaks

\def\strut{\vphantom{\Large(} }

\def\tat{t\^ete-a-t\^ete}
\def\init{\operatorname{init}}
\def\ini{\operatorname{in}}
\def\Ker{\operatorname{Ker}}
\def\supp{\operatorname{supp}}
\def\Sec{\operatorname{Sec}}

\def\Sagbi{{\textsc{Sagbi}}}

\begin{document}


\title{Sagbi bases, defining ideals and algebras of minors}

\author{Winfried Bruns}
\address{Universit\"at Osnabr\"uck, Institut f\"ur Mathematik, 49069 Osnabr\"uck, Germany}
\email{wbruns@uos.de}
\author{Aldo Conca}
\address{Dipartimento di Matematica, Dipartimento di Eccellenza 2023-2027, Universit\'a di Genova, Genova, Italy}
\email{aldo.conca@unige.it}
\author{Francesca Lembo}
\address{Dipartimento di Matematica, Dipartimento di Eccellenza 2023-2027, Universit\'a di Genova, Genova, Italy}
\email{francescalembo@outlook.com}


\begin{abstract}
This paper extends the article of Bruns and Conca on \Sagbi{} bases and their computation (J. Symb. Comput. 120  (2024)) in two directions. (i) We describe the extension of the Singular library sagbiNormaliz.sing to the computation of defining ideals of subalgebras of polynomial rings. (ii) We give a complete classification of the algebras of minors for which the generating set is a \Sagbi{} basis with respect to a suitable monomial order and we identify universal  \Sagbi{} basis in three cases. The investigation is illustrated by several examples.
\end{abstract}

\thanks{A.C. is supported by PRIN~2020355B8Y ``Squarefree Gr\"obner degenerations, special varieties and related topics,'' by MIUR Excellence Department Project CUP~D33C23001110001, and by INdAM-GNSAGA}

\maketitle
 
 \section{Introduction}

Let $R$ be a polynomial ring over a field $K$, endowed with a monomial (or term) order and $A$ be a $K$-subalgebra of $R$. The $K$-vector space spanned by the initial (or leading) monomials of the elements of $A$  is a $K$-subalgebra  $\ini(A)$ of $R$. A \Sagbi{} basis is a subset of $A$ whose initial monomials generate $\ini(A)$ as a $K$-algebra. There is a clear analogy to Gr\"obner bases of ideals $I\subset R$: a subset $G$ of $I$ is a Gr\"obner basis if the initial monomials of the polynomials $f\in G$ generate (as an ideal) the initial ideal of $I$. The acronym \Sagbi{}, introduced by Robbiano and Sweedler \cite{RobSwe},  stands for ``subalgebra analog to Gr\"obner bases of ideals''. Section \ref{basics} gives a compact introduction to \Sagbi{} bases and their computation.

In \cite{BCSagbi} Bruns and Conca have described the implementation of the \Sagbi{} algorithm in the Singular \cite{Sing} library sagbiNormaliz.lib which uses Singular as the environment for polynomial computations and Normaliz \cite{Nmz} for the combinatorics. The first major goal of this note is the extension of the library by functions that, in addition to a \Sagbi{} basis, compute a defining ideal of the subalgebra $A$ in terms of the given generators (see Section \ref{DefId}). This extension is possible since \Sagbi{} bases are characterized by the liftability of the binomial ideal defining the initial algebra $ \ini(A)$ to a defining ideal of $A$. Section \ref{Compute} compares computation times of the approach via \Sagbi{} bases to the classical elimination by CoCoA-5 \cite{CoCoA} and Singular. No doubt, in general elimination is faster, but there are interesting cases in which the \Sagbi{} approach is competitive or better.

Given positive integers $m,n$ with $m\leq n$ let $X_{m\times n}$ be a $m\times n$  matrix of variables over a field $K$ and $R=K[X_{m\times n}]$ be the polynomial ring  generated by its  entries.  The algebras, whose \Sagbi{} bases are most interesting to the authors, are the subalgebras $A_t(m,n)$ of $R$ generated by the $t$-minors of $X_{m\times n}$. The $t$-minors are the determinants of $t\times t$ submatrices of $X$. In the case $t = m$, the algebra  $A_m(m,n)$ is also denoted by $G(m,n)$ and is  the homogeneous coordinate ring of the classical Grassmannian, i.e., the variety of $m$-dimensional $K$-subspaces of the vector space $K^n$.

In Section \ref{tools} we recall notions and tools introduced by Sturmfels and Zelevinsky \cite{SZ}, \cite{StuGrPol} for the investigation of the collections of initial monomials of the $m$-minors in $G(m,n)$ under varying monomials orders. While the maximal minors do not form a \Sagbi{}  basis of $G(m,n)$ in general, the collections of initial monomials always generates a subalgebra of the same dimension as $G(m,n)$. This is not only of theoretical interest, but has computational consequences that we will use in later sections. The main tool for this result is the theory of Cartwright-Sturmfels ideals developed by Conca, De Negri and Gorla \cite{CDG1}--\cite{CDG6}.

It is known by a theorem of Sturmfels \cite{StuGrPol} that the $m$-minors form a \Sagbi{}  basis of $G(m,n)$ with respect to a diagonal monomial order, i.e., a monomial  order for which the initial monomials of the  minors  are the product of the elements in their  main diagonal. By a theorem of Bruns and Conca (for example, see \cite[Sect. 6.4]{BCRV} or \cite[Thm. 3.11]{BC}) if the characteristic of $K$ is not too small, $A_t(m,n)$ always has a finite \Sagbi{} basis for a diagonal order. But, apart from the trivial case $t=1$ and the case of Grassmannians, the $t$-minors do not form a \Sagbi{} basis for a diagonal order. This raises the question whether there exists a monomial order on  $K[X_{m\times n}]$ for which the $t$-minors are a \Sagbi{} basis. The main result of Section \ref{AreThey?} is Theorem \ref{AtSagbi}: this holds if and only if $t = m$, or  $t+1=m=n$ or $t = 1$. In the key case $t+1 = m = n$  the selection of initial monomials that makes the $(m-1)$-minors a \Sagbi{} basis  is unique up to the permutations of rows and columns (Theorem \ref{unique}). In turn this follows from the Birkhoff--von Neumann theorem, applied to the Newton polytope of the product of $(m-1)$-minors of $X_{m\times m}$. We can then show that for larger formats the $t$-minors cannot be a \Sagbi{} basis.

We illustrate the technique of the Newton polytope by several concrete examples, and demonstrate the use of Hilbert series in (dis)proving that certain sets of minors or products of minors are a \Sagbi{} basis. This allows to find a universal \Sagbi{} basis for $m=n=3$ and $t=2$ (Theorem \ref{univSagbi3x3}), for $G(3,6)$ (Theorem \ref{univSagbi3x6}) and $G(3,7)$ (Theorem \ref{univSagbi3x7}).   Hilbert series computations for the algebras generated by initial monomials of the $3$-minors in $G(3,8)$ and $G(3,9)$ show that universal \Sagbi{} bases for them must be very complicated. 

As a potential  extension in another direction, in Remarks \ref{R=G}--\ref{Rees_2} we discuss the problem of finding universal \Sagbi{} bases for the Rees algebra of the ideal  of $3$-minors of  $X_{3\times 6}$ and  of the ideal  of $3$-minors of  $X_{3\times 7}$. The first case seems to be treatable while the second appears to be  very complicated.  

We point out that   \Sagbi{} bases and initial algebras of $G(m,n)$ and, in particular, of $G(3,n)$ have been studied by several authors,  see \cite{BMC, CM1, CHM, CM2, MS, M}. 

The results in Sections \ref{AreThey?} and \ref{univSag}  generalize statements contained in  the third author's master thesis \cite{Lembo} written under the supervision of the second author. 

We would like to thank Barbara Betti and the anonymous reviewers for their insightful comments on an earlier version of the paper, which helped us improve the quality of the exposition.  

\section{Basics of \Sagbi {} bases}
\label{basics}

This section is  a slightly modified version  of \cite[Sect. 2]{BCSagbi}. It is included to keep this paper as self-contained as possible.

The reader finds a compact discussion of \Sagbi {} bases in \cite[Sect. 1.3]{BCRV}. We use the notation developed there. Kreuzer and Robbiano \cite[Sect. 6.6]{KrRo} give a more extensive introduction;  see also Ene and Herzog \cite{EneHerz} and Sturmfels \cite{StuGrPol}. \Sagbi {} bases were introduced independently by Robbiano and Sweedler \cite{RobSwe} and Kapur and Madlener \cite{KapMad}. The acronym \Sagbi {} stands for ``subalgebra analog to Gr\"obner bases of ideals'' \cite{RobSwe}.     Some authors have adopted recently a new terminology, Khovanskii bases, for a notion that generalize that \Sagbi {} bases but  in this paper we keep the traditional name.

Let $A \subset R=K[X_1,\dots,X_n]$ be a $K$-subalgebra and $\cF$ a (not necessarily finite) family of polynomials belonging to $A$. We assume that $R$ is endowed with a monomial (or term) order $<$. In the following  we often employ  a  simplified terminology using \emph{subalgebra}  for $K$-subalgebra   and \emph{order} for monomial order. One calls $\cF$ a \emph{\Sagbi {} basis} of $A$ if the initial monomials $\ini(f)$, $f\in \cF$, generate the initial algebra $\ini(A)$. A \Sagbi {} basis is automatically a system of generators of $A$. If $\cF$ is finite, then $A$ and $K[\ini(\cF)]$ are connected by a flat deformation (Conca, Herzog and Valla \cite{CHV}), and this allows the transfer of homological properties and numerical invariants  from the toric algebra $K[\ini(\cF)]$ to $A$. Chapter 6 of \cite{BCRV} exploits this approach for the investigation of algebras generated by minors. 

\Sagbi{} bases need not be finite, but can always be chosen countable. Therefore one must allow that $\cF= (f_u)_{u\in N}$ with $N=\{1,\dots,p\}$ with $p\in \NN$ or $N=\NN$. We will always assume that the members of $\cF$ are monic. This is evidently no essential restriction of generality as long as the base ring $K$ is a field. However, in the bookkeeping underlying the computation of the defining ideal the division of a polynomial by its initial coefficient must of course be registered.

For us the following simple lemma is an important tool in the computation of \Sagbi{} bases. For an $\NN$-graded $K$-vector space $V$ one defines the \emph{Hilbert function} of $V$ by
\begin{equation}
H(V,k) = \dim_K V_k, \qquad k\in\NN,\label{Hilb}
\end{equation}
where $V_k$ is the subspace of degree $k$ elements of $V$.

\begin{lemma}\label{HilbLemma}
Let $R$ be a polynomial ring, endowed with a monomial order and the standard grading. Let  $A$ be a finitely generated graded subalgebra. Furthermore let $\cF$ be a family of polynomials in $A$ and $B=K[\ini(\cF)]$ the subalgebra generated by the monomials $\ini(f)$, $f\in \cF$. Then the following hold:
\begin{enumerate}
	\item $H(B,k) \le H(A,k)$ for all $k$.
	\item $\cF$ is a \Sagbi {} basis of $A$ if and only if $H(B,k) = H(A,k)$ for all $k\in\NN$.
\end{enumerate}
\end{lemma}

\begin{proof}
For a graded subspace $V$ of the polynomial ring one has $H(V,k)=H(\ini(V)),k)$ for all $k$ \cite[1.4.3]{BCRV}. Furthermore there are inclusions
$$
B_k \subset \ini(A)_k,\qquad k\in\NN,
$$
and equality holds if and only if $H(B,k)=H(\ini(A),k)$ for all $k$. Together with Equation \eqref{Hilb} this proves the lemma.
\end{proof}

To present $A$ as a residue class ring of a polynomial ring, we choose $P=K[Y_u: u\in N]$ and define a surjection
\begin{equation}
\phi: P \to A, \qquad \phi(Y_u) = f_u, \ u\in N.\label{phi}
\end{equation}
The $K$-algebra  $K[\ini(\cF)]$  is a homomorphic image of $P$, as well, namely by the surjection
$$
\psi: P     \to K[\ini(\cF)], \qquad \psi(Y_u) = \ini(f_u), \ u\in N.
$$
The kernel of $\psi$ is generated by a set of binomials. In the terminology of \cite{RobSwe}, a binomial in $\Ker\psi$ is called a \emph{\tat}.

A monomial in $P$ is given by an exponent vector  $e=(e_u)_{u\in N}$ of natural numbers $e_u$ of which all but finitely many are $0$. We set $Y^{e} = \prod_{u\in N} Y_u^{e_u}$ and 
$$\cF^e=\psi(Y^{e})=\prod_{u\in N} f_u^{e_u}.$$
Let $F\in P$ be a polynomial, given as a finite $K$-linear combination of monomials, $F=\sum_i a_iY^{e_i}$ where the $e_i$ are exponent vectors, and $a_i\neq 0$ for all indices involved. We set
\begin{align*}
\ini_\phi(F)&=\max_i \  \ini(\phi(Y^{e_i})), \\
\init_\phi(F)&=\sum_{\ini(\phi(Y^{e_i}))=\ini_\phi(F) } a_iY^{e_i}.
\end{align*}
Note that in the definition of $\ini_\phi(F)$ the maximum is taken over the initial monomials with respect to the monomial order on $R$ so that it is a monomial in $R$. In contrast, $\init_\phi(F)$ is a polynomial in $P$. Since $\phi(\init_\phi(F))$ can be $0$, in general $\ini_\phi(F) \neq \ini(\phi(F))$, and this cancellation of initials is the crucial  aspect of \Sagbi {} computation.  One says that a polynomial $F\in\Ker\phi$ \emph{lifts} a polynomial $H\in \Ker\psi$ if $\init_\phi(F)=\init_\phi(H)$.

We can now formulate the \emph{\Sagbi {} criterion} (see \cite[1.3.14]{BCRV}).

\begin{theorem}\label{criterion}
With the notation introduced, let $\cB$ be a set of binomials generating $\Ker\psi$. Then the following are equivalent:
\begin{enumerate}
\item $\cF$ is a \Sagbi {} basis of $A$;
\item every binomial $\beta\in \cB$ can be lifted to a polynomial $F_\beta\in \Ker\phi$.
\end{enumerate}
If $\cF$ is a \Sagbi {} basis of $A$, then the polynomials $F_\beta$ generate $\Ker \phi$. 
\end{theorem}

The Buchberger algorithm for a Gr\"obner basis of an ideal $I$ starts from a system of generators $G$ of $I$. Then one applies two steps, namely (i) the computation of the $S$-polynomials $S(g_1,g_2)$, $g_1, g_2\in G$, and (ii) their reductions modulo $G$. The nonzero reductions are then added to $G$, and the next round of $S$-polynomials of the augmented $G$ and their reductions is run. This produces an increasing sequence of initial ideals $\ini(G)$. Because of Noetherianity the process stops after finitely many rounds with a Gr\"obner basis of $I$. The reduction of an $S$-polynomial $S(g_1,g_2)$ to $0$ is equivalent to the liftability of the  ``divided Koszul syzygy'' of $\ini(g_1)$ and $\ini(g_2)$ to a syzygy of the polynomials $g\in G$.

The computation of \Sagbi {} bases follows the same pattern. There are however two main differences: an analog of the divided Koszul syzygies does not exist, and ascending chains of monomial subalgebras of $R$ need not stabilize. For \Sagbi {} bases one must therefore compute a binomial system of generators of $\Ker\Psi$, and one cannot expect the algorithm to stop. The analog of reduction is called subduction (we again follow \cite{RobSwe}) .

\begin{definition}\label{subduction}
Let $g\in R$. Then $r\in R$ is a \emph{\index{subduction}subduction of $g$ modulo $\cF$} if there exist monomials $\cF^{e_1},\dots,\cF^{e_m}$ and non-zero coefficients $a_i\in K$ such that the following hold:
\begin{enumerate}
\item $g=a_1\cF^{e_1}+\dots+a_m\cF^{e_m}+r$;
\item $\ini(\cF^{e_i})\le \ini(g)$ for $i=1,\dots.m$;
\item no monomial $\mu\in\supp(r)$ is of type $\ini(\cF^e)$.
\end{enumerate}
\end{definition}
The process that computes a subduction of $g$ modulo $\cF$ is also called subduction. Here $\supp(r)$ is the set of monomials of $R$ appearing in $r$ with a nonzero coefficient. In the computation of \Sagbi {} bases,  including defining ideals,  one can replace (3) by the weaker condition
\begin{enumerate}
\item[($3'$)]  $\ini(r)$ is not of type $\ini(\cF^e)$.
\end{enumerate}
There is an obvious algorithm that produces a \emph{subduction remainder} $r$ in $(3')$: if $\ini(g)=\ini(\cF^e)$, we   replace $g$ by $g - a\phi(\cF^e)$ where $a$ is the leading coefficient of $g$, and iterate this \emph{subduction step} as long as possible. The algorithm stops since the sequence of initial monomials is descending and descending sequences in a monomial order are finite. Once $(3')$ is reached, one applies subduction steps iteratively to the remaining monomials to achieve the ``tail subduction'' asked for by (3). 

The crucial observation for the computation of the defining ideal and for the proof of correctness of the \Sagbi{} algorithm below is the following proposition.

\begin{proposition}\label{lift}
With the notation above, suppose that  $\phi(\beta)$ subducts   to $0$ for a \tat{} $\beta\in \Ker\psi$, $\phi(\beta) = a_1\cF^{e_1}+\dots+a_m\cF^{e_m}$. Then
$$
\beta - (a_1Y^{e_1}+\dots+a_mY^{e_m}) \in \Ker \phi 
$$
lifts $\beta$.
\end{proposition}

It is enough to observe $\ini(\cF^{e_i})\le\ini(\phi(\beta)) < \ini_\phi(\beta)$ for all $i$.

The algorithm (\emph{\Sagbi {}}) starts from the finite family $\cF_0$ generating the subalgebra $A\subset R$. Then one proceeds as follows:
\begin{enumerate}
\item[(1)] Set $i=0$.
\item[(2)] Set $\cF'=\emptyset$ and compute a binomial system of generators $\cB_i$ of the kernel of 
$\psi_i: P_i \to K[\ini(\cF_i)]$, $P_i=K[Y_F:F\in \cF_i]$, $\psi_i(Y_F)=\ini(F)$.
\item[(3)] For all $\beta\in \cB_i$ compute the subduction $r$ of $\phi_i(\beta)$ modulo $\cF_i$, $\phi_i$ given by the substitution $Y_F\mapsto F$, $F\in\cF_i$.  If $r\neq 0$, make $r$ monic and add it  to $\cF'$.
\item[(4)] If $\cF'=\emptyset$, set $\cF_j=\cF_i$, $P_j=P_i$, $\cB_j=\cB_i$ for all $j\ge i$ and stop.
\item[(5)] Otherwise set $\cF_{i+1}=\cF_i\cup \cF'$, $i=i+1$ and go to (2).
\end{enumerate}
It follows from Theorem \ref{criterion} and Proposition \ref{lift} that $\cF=\bigcup_{i=0}^\infty\cF_i$ is a \Sagbi {} basis of $A$. It is not hard to see that the algorithm stops after finitely many steps if $A$ has a finite \Sagbi {} basis (with respect to the given monomial order). 

To complete the notation we define the polynomial ring $P_\infty = \bigcup_{i=0}^\infty P_i$. It is endowed with surjective ring homomorphisms $\phi_\infty: P_\infty  \to A$, $\phi_\infty(Y_F) = \phi_i(Y_F) = F$ for $F\in \cF_i$, and the analogously constructed $\psi_\infty: P_\infty \to \ini(A)$, $\psi_\infty(Y_F) = \psi_i(Y_F) = \ini(F)$.

\section{Computing the defining ideal}
\label{DefId}

The algebra $A$ is given by the family of generators $\cF_0$. The polynomial ring $P_0 = K[Y_F: F\in \cF_0]$ is mapped surjectively onto $A$ by the $K$-algebra map $\pi: P_0\to A$, $\pi(Y_F)= F$, $F\in \cF_0$.  We write $\pi$ instead of $\phi_0$ for better readability. Computing the defining ideal means to compute a (finite) system of polynomials in $P_0$ generating the kernel of $\pi$. 

Theorem \ref{criterion} and Proposition \ref{lift} show that we know a defining ideal of $A$: $\Ker\phi_\infty$ is generated by the polynomials that represent the subduction to $0$ of the \tat{}s. But this defining ideal lives in the  ``wrong'' polynomial ring, namely in $P_\infty$ instead of $P_0$. The ``bookkeeping'' that we need for computing $\Ker \pi$ from $\Ker\phi_\infty$ is a retraction  $\rho: P_\infty \to P_0$ of the natural embedding $P_0\hookrightarrow P_\infty$ such that $\phi_\infty = \pi \circ \rho$. Evidently, $\Ker\pi =\rho(\Ker\phi_\infty)$ then.

We construct $\rho$ and $\Ker\pi$ inductively. Suppose we have reached step (3) in the \Sagbi{} algorithm, and let  $\phi_i(\beta) = a_1\cF^{e_1}+\dots+a_m\cF^{e_m} +r $ be a subduction of $\phi_i(\beta)$. There are two cases:
\begin{enumerate}
\item If $r\neq 0$, then we register
$$
\rho(r) = \rho(\beta - (a_1Y^{e_1}+\dots+a_mY^{e_m})).
$$
This makes sense since $\rho$ has already been computed on $P_i$. To make $r$ monic, divide both sides by the initial coefficient of $r$.

For the indeterminate $Y_r$ introduced in $P_{i+1}$ we set $\rho(Y_r) = \rho(r)$.  After all nonzero subduction remainders have been done, the retraction $\rho$ has been extended to $P_{i+1}$.

\item If $r = 0$, then  $\rho(\beta - (a_1Y^{e_1}+\dots+a_mY^{e_m}))$ is added as a generator to $\Ker\pi$.
\end{enumerate}

\begin{remark}
(a) Instead of using a retraction  from $P_i$, $i\ge0$ to $P_0$ one could consider a chain of retractions $\rho_{i+1}: P_{i+1} \to P_i$, and go down the chain at the end. However, it seemed easier to implement the ``all at once'' method above.

(b) In addition to the defining ideal the library also returns the retraction $\rho$, in other words, expressions of the \Sagbi{} basis elements in terms of the variables in $P_0$. These are often much shorter than the representations of the elements as polynomials in $R$. 
\end{remark} 

For the computation of the \Sagbi{} basis without the defining ideal the Singular library sagbiNormaliz.lib has $3$ variants of the \Sagbi{} algorithm:
\begin{enumerate}
\item[(Gen)] uses the \Sagbi{} algorithm as given above and can be applied to all subalgebras of polynomial rings.
\item[(Deg)] can be applied to graded subalgebras of standard graded polynomial rings. In computing the \Sagbi{} basis  it proceeds degree by degree.
\item[(Hilb)]can be applied in the same cases as (Deg), provided the Hilbert series is known. It proceeds by degrees, using that the next ``critical'' degree can be read from the Hilbert function. Moreover, it avoids subductions as soon as it is clear that, roughly speaking, as many ``new'' elements of the \Sagbi{} basis have been found as predicted by  the difference of Hilbert functions of $A$ and the subalgebra of $A$ generated by the subset of the \Sagbi{} basis computed so far.
\end{enumerate}
See \cite[Sect. 4]{BCSagbi} for more precise information on these variants. 

Since (Hilb) avoids complete subductions, it cannot be used for the computation of the defining ideal. Therefore only augmented versions of (Gen) and (Deg) have been provided for it. For the augmented versions, the same rule of thumb applies as for the pure \Sagbi{} computations: whenever possible, use (Deg).

We use the \emph{normalized degree} on the subalgebra $A$: it is obtained from the standard degree of the surrounding polynomial ring by division by the $\gcd$ of the standard degrees of the generators. 

For the combinatorial computations involving monomial algebras and binomial ideals the library falls back on Normaliz. These computations do not change when the defining ideal is asked for. This is the first reason why computing the defining ideal extends the computation only mildly in most cases. Another reason is that the computation of the retraction $\rho$ uses indeterminates instead of the potentially complicated original generators of the subalgebra~$A$. 

\section{Computational data}
\label{Compute} 

The subalgebras to which we have applied our algorithm are generated by minors of a matrix of inderterminates or closely related to them. Using the shortcut
$$
[m]= \{1,\dots m\},
$$
we let
$$
K[X] =K[X_{m\times n}] = K\bigl[X_{ij}: i\in [m],\ j\in [n]\bigr]
$$
denote the polynomial ring in $mn$ indeterminates, arranged in an $m\times n$ matrix. We assume that $m\le n$, unless stated otherwise. The standard monomial orders on such a polynomial ring are the \emph{diagonal} ones: for such an order, the initial monomial of any  \emph{$k$-minor}  is the product of the indeterminates in the main diagonal of the submatrix.  In this section, unless stated otherwise, a diagonal monomial order is used. Our standard reference is Bruns, Conca, Raicu and Varbaro  \cite{BCRV}. 

We have tested the computation of defining ideals for algebras of the following types.
\begin{description}
\item[$G(m,n)$] homogeneous coordinate ring of the Grassmannian of $m$-spaces in $K^n$ with respect to the Plücker embedding. It is subalgebra of the polynomial ring $K[X_{m\times n}]$, generated by the $m$-minors. It is well-known that the $m$-minors not only generate $G(m,n)$, but even are a \Sagbi{} basis \cite[Thm. 6.2.1]{BCRV}. The defining ideal is generated by the degree $2$ \emph{Plücker relations} \cite[Sect. 3.2]{BCRV}. These properties are independent of the characteristic of the base field. In the computations we use characteristic $0$.

\item[$G(m,n)_\ell$] the same as $G(m,n)$, but with a non-diagonal lexicographic order.

\item[$A_t(m,n)_c$] is the subalgebra of $K[X]$, $\chara K =c$, generated by the $t$-minors. If $c= 0$ or $c > \min(k,m-t,n-t)$, it has a finite \Sagbi{} basis of which the $t$-minors form a strict subset in general \cite[Sect. 6.4]{BCRV}. The defining ideal is unknown in general, but for characteristic $0$ conjectured to be generated in degree $2$ and $3$, see Bruns, Conca and Varbaro \cite{BCV} or \cite[Sect. 6.5]{BCRV}. For $t= 2$ this conjecture has been proved by Huang, Perlman, Polini, Raicu, and Sammartano \cite{HPPRS}.

\item[$R_{r}(m,n)$] Let $Y$ be an $m\times r$ matrix of indeterminates and $Z$ an $r\times n$-matrix. The subalgebra of $K[Y,Z]$ generated by the entries of the product matrix $YZ$ is isomorphic to the residue class ring of $K[X_{m\times n}]$ by the ideal $I_{r+1}$ generated by the $(r+1)$-minors. Computing its \Sagbi{} basis amounts to computing a toric deformation of the determinantal variety of rank $r$ matrices of format $m\times n$. See \cite[Sect. 6.7]{BCRV}. From the description given it is clear that the defining ideal of this algebra is $I_{r+1}$.

\item[$\Sec_k(m,n)$] denotes the homogeneous coordinate ring of the $k$th secant variety of $G(m,n)$. See \cite{AOP}.

\end{description}

Table \ref{Bench_1} lists computation times, cardinalities, and degrees of \Sagbi{} bases and defining relations of several examples. We compare the computation times of sagbiNormaliz.lib to classical elimination by Singular and CoCoA5. 
  
\begin{table}[hbt]
\begin{tabular}{crrrrrrrr}
\midrule[1.2pt]
\strut        & \multicolumn{3}{c}{norm deg}& & &\multicolumn{3}{c}{times in minutes}\\
\cline{2-4}\cline{7-9}
\strut ecample& bound& \Sagbi {} & comp & \#\Sagbi {} & \#rel& (Deg) & Sing & CoCoA \\
\midrule[1.2pt]
\strut $A_2(4,4)_0$ & 10&    3 & 7 & 89 & 240 &0:27.2 & 0:10.2 & 0:47.3 \\
\hline
\strut $A_2(4,4)_0$ &  3&    3 & 3 & 89 & 240 &0:10.7 & -- & -- \\
\hline
\strut $A_2(4,4)_2$ & 15&    6 & 13 &130& 205 & 4:09.6 & 0:11.0 & 0:49.3  \\
\hline
\strut $R_2(4,4)$ & 15 &  4 & 5 & 52 & 16 & 0:01.2 & 0:00.9 & 0:00.5\\
\hline
\strut $G(3,7)$ & 15 & 1 & 3 & 35 & 140 & 0:06.1 & 1:16.2 & 0:01.1\\
\hline
\strut $G(3,8)$ & 15 & 1 & 3 & 56 & 420 & 0:02.2 & 432:46.0  & 0:08.4\\
\hline
\strut $G(3,9)$ & 15 & 1 &3 & 84 & 1050 & 0:10.2 & $> 1$ d & 1:01.1\\
\hline
\strut $G(3,11)$ & 15 & 1 &3& 165 & 4620 & 24:05.7 & $> 1$ d & 40:57.1 \\
\hline 
\strut $G(3,7)_\ell$ & 15 & 2 & 4 & 37 & 140 & 0:00.9 & 01:15.4 & -- \\
\hline 
\strut $G(3,9)_\ell$ & 15 & 3 & 6 & 101 & 1050 & 0:27.8 & $> 1$ d & -- \\
\hline
\strut $\Sec_2(3,7)$ & 3 & 3 & 3& 284 & 28 & 247:26.0 &  $> 1$ d &  $> 1$ d\\
\hline
\end{tabular}
\vspace*{1ex}
\caption{sagbiNormaliz.lib vs. elimination by Singular and CoCoA}	
\label{Bench_1}
\end{table}

Table \ref{Bench_1} contains only graded algebras. Consequently we have used the degree-by-degree variant of the algorithm.  

In the table the column ``bound'' is rather irrelevant. The degree-by-degree variant asks for it. In the cases where we could compute a complete \Sagbi{} basis it was chosen large enough to allow the complete computation.  In the case of  $\Sec_2(3,7)$ the degree bound was set to $3$ because a previous computation with t a larger degree bound suggested that  no new relations came up after that degree.

for which a previous computation with a larger bound suggested the relations in degree $3$.

The column ``\Sagbi'' lists the maximal degree of the elements of the \Sagbi{} basis (as far as computed). The column ``comp'' shows the degree to which the computation runs until the completion of the \Sagbi{} basis is recognized or the degree bound is reached. These numbers seem to be too large by $1$ in cases like $G(m,n) $ where it is known a priori that the \Sagbi{} basis lives in degree $1$. However, the implementation goes degree by degree until it reaches a point at which  the  \tat{} has no binomials anymore.

The next two columns list the number of elements of the \Sagbi{} basis and relations, as far as computed.

All computations were done on a PC with a AMD Ryzen 9000 CPU.
\begin{remark}

(a) It is not surprising that in many cases classical elimination beats the approach via \Sagbi{} bases. The former always terminates, whereas the latter yields a definite result only if the \Sagbi{} basis is finite or some extra information can be used, in particular a degree bound. An example is $A_2(4,4)_0$  for which the generation in degrees $2$ and $3$ has been proved in \cite{HPPRS}, and an explicit list of generators is asked for.

Also the defining ideal of $A_2(4,4)_2$ is generated in degrees $2$ and $3$. The \Sagbi{} basis is however much larger than the one in characteristic $0$.

(b) $A_2(4,4)_2$ is a case in which  the computation of the defining ideal essentially doubles the computation time of the \Sagbi{} basis. This is almost entirely due to computation time needed for the final minimization of the generating system of relations.

(c) The main point in trying $G(m,n)$ is finding an explicit system of generators for the defining ideal without much preparation. The computation via \Sagbi{} basis does indeed offer a fast approach. We have no explanation for the explosion of computation times for classical elimination. An essential point may be that for the \Sagbi{} basis a Gröbner basis computation is only necessary for a single binomial ideal (in this case) followed by rather fast computation of subductions in degree $2$.

(d) $\Sec_2(3,7)$ is a very interesting case. Whereas classical elimination did not give any information on the defining ideal, not even with a degree bound, in reasonable time, the \Sagbi{} approach succeeded in a few hours by listing $27$ degree $2$ relations. One cannot conclude from the computation that this system of generators is complete. However, the ideal generated by these $27$ polynomials defines a residue class ring of the right dimension and the numerator polynomial of the Hilbert series has positive entries.  
\end{remark}


\section{Coherent matchings, Newton polytopes and computational tools}\label{tools}

For the computational proofs in the next sections we recall and extend some terminology and facts  from  \cite{StuGrPol}. For a  polynomial $F\in K[X_1,\dots, X_n]$ we denote by $\New(F)$ the Newton polytope of $F$, i.e., the convex polytope  spanned  by the exponent vectors of the monomials in the support of $F$. 

For a finite list  $\cF=\{F_1,\dots, F_s\}$ of nonzero polynomials we let $\New(\cF)$ denote the Newton polytope $\New(F_1\cdots F_s)$ of  the product of the elements in $\cF$. 
For such a list $\cF$ a \emph{matching} $T$ is a list of monomials $\{X^{a_1}, \dots, X^{a_s} \}$ such that   $X^{a_i}$ is in  the support of  $F_i$ for $i=1,\dots, s$. A matching is \emph{coherent} if it consists of the initial monomials with respect to some monomial  order or, equivalently, of some general enough positive  integral weight vector.  (For example, see \cite[Lemma 1.5.6]{BCRV}.) 

A vertex of $\New(\cF)$  can be given by specifying a sufficiently general weight vector $w\in \ZZ^n$. The vertex $v$ of $\New(\cF)$ is the unique point of $\New(\cF)$ that maximizes $w$. Without further hypothesis, it may not be possible to choose $w$ positive. But if the elements in  $\cF$ are homogeneous (under an arbitrary positive grading), then $w$ can be chosen as a positive vector. Indeed  the vector $g$ formed by the degrees of the variables is constant on every homogeneous polynomial, and $v$ also maximizes $w + mg$ for $m \gg 0$, which has positive entries.   An important observation here is that $\New(\cF)$ is the Minkowski sum of the  polytopes $\New(F_i)$ with $i\in [s]$ and hence each vertex $v$ of  $\New(\cF)$  is the sum $v_1+\dots v_s$, in a unique way,  of vertices $v_i$ of  $\New(F_i)$.   Summing up we have

\begin{proposition} 
Assume the polynomials in  $\cF=\{F_1,\dots, F_s\}$  are homogeneous. Then   there is a one-to-one correspondence between the coherent matchings of the list  $\cF$ and the vertices of the Newton polytope $\New(\cF)$ of the product $F_1\cdots F_s$. 
\end{proposition}

\begin{remark}
With the notation introduced, let $A$ be the subalgebra of $K[X_1,\dots,X_n$ ] generated by $\cF$. We say that a coherent matching $T$ of $\cF$ or, equivalently, the corresponding vertex of the Newton polytope  $\New(\cF)$ is \Sagbi{} if $\cF$ is a \Sagbi{} basis of $A$ with respect to a monomial order that selects $T$ from the elements of $\cF$. The definition makes sense since this property does not depend on the chosen order. In fact, $K[T]$ is the initial algebra of $A$ with respect to at least one monomial order, and is contained in $\ini_<(A)$ for any monomial order $<$ that picks $T$ from $\cF$. But then $K[T] = \ini_<(A)$, as follows from \cite[Prop. 1.4.1(e)]{BCRV}. 
\end{remark}

The observations above allow us to make statements about the \Sagbi{} property of $\cF$  with respect to all monomial orders: their number is infinite, but the Newton polytope has only finitely many vertices, or equivalently, the number of coherent matchings is finite. Let us list the computational tools that we will employ in the following. In Section \ref{AreThey?} these tools will only be used to provide illustrative examples, but in Section \ref{univSag} the proofs will be based on them. 
\begin{enumerate}
\item We have used CoCoA \cite{CoCoA} for computations based on random weight vectors, especially in connection with the computation of Hilbert functions. 

\item Normaliz \cite{Nmz} has been used to confirm and extend the results of (1) on vertices of Newton polytopes.

\item We have extended a C++ program based on libnormaliz and CoCoALib that we already used in \cite{BCSagbi} for the computation of all coherent matchings of the set of $3$-minors in $G(3,n)$ and deciding their \Sagbi{} property. In the following we call this tool \emph{SagbiGrass}.
 \end{enumerate}
 
\begin{remark}\label{correction}
(a) A new check has shown that \cite[Rem. 9(b)]{BCSagbi} is false: already for $G(3,6)$ there exist coherent matchings that are neither lex nor revlex compatible. 

(b) \cite[Thm.\ 7(2)]{BCSagbi} can be improved: the algebra generated by the initial monomials of a coherent matching of $G(3,n)$ is always normal for $n  \le 8$. By Hochster's \cite[Thm.6.3.5]{BH}theorem it follows that all these algebras are Cohen--Macaulay and therefore have positive coefficient vectors in the numerators of their Hilbert series. As already stated in  \cite[Thm.\ 7(3)]{BCSagbi} normality sometimes fails for $n = 9$.
\end{remark}


\section{Dimension of coherent matchings for $G(m,n)$}
\label{dimension_ok}
Denote by $\cM_m$ the set of $m$-minors of  $X_{m\times n}$.
Let $T$ be coherent matching  of $\cM_m\subset K[X_{m\times n}]$ and let $K[T]$  be  the $K$-algebra it generates.   Then by Lemma \ref{HilbLemma} we have $H(K[T],i)\leq H(G(m,n),i)$ for all $i\in\NN$  and  $T$ is \Sagbi{} if and only if equality holds for all $i\in \NN$. Hence the difference $H(G(m,n),i)-H(K[T],i)$ can be taken as a measure of the failure of $T$ being \Sagbi{}.  The goal of this section is to prove that the degree of growth of the $H(G(m,n),i)$ and $H(K[T],i)$ is always the same, that is: 

\begin{theorem}  \label{dimok}
 For every coherent matching $T$ of $\cM_m\subset K[X_{m\times n}]$ one has $$\dim K[T]=\dim G(m,n).$$
 \end{theorem}

 In the proof we will use Cartwright-Sturmfels ideals and their   properties. This theory has been developed in the series of papers \cite{CDG1,CDG2,CDG3,CDG4,CDG5,CDG6} and has  roots and applications in \cite{Bo,CS,C, CW}.   
 
 \begin{proof}  The ideal $I_m(X)\subset K[X_{m\times n}]$ generated by $\cM_m$ is a 
 Cartwright-Sturmfels ideal  with respect  to the row $\ZZ^m$-grading, i.e.~$\deg X_{ij}=e_i\in \ZZ^m$, see \cite[Thm.1]{CDG3}. Its  $\ZZ^m$-graded generic initial ideal,  with respect to any monomial order $<$ that satisfies  $X_{ij}>X_{ik}$ for all $i$ and for all $j<k$,  is generated by the elements in the set 
 $$B_{m,n}=\left \{ \prod_{i=1}^m X_{ij_i} \  : \  j_1+j_2+\cdots +j_m \leq n \right\},$$
 see \cite[Sect.4]{CDG1} or \cite[Thm.1]{CDG3}. 
 Let $J$  be the ideal of $K[X_{m\times n}]$ generated by $T$. Its degree $m$ component $J_m$    is  the $K$-vector space generated by $T$.  Since $\cM_m$ is a universal Gr\"obner basis \cite{SZ, BZ, CDG1}, the ideal $J$ is an initial ideal of the ideal $I_m(X)$.
It  follows from \cite[2.4, 2.5]{CDG2} that $J$  is a Cartwright-Sturmfels ideal  as well and that its $\ZZ^m$-graded generic initial ideal  equals that of $I_m(X)$.    Let $\phi$ be a generic $\ZZ^m$-graded $K$-algebra automorphism of $K[X_{m\times n}]$. The fact that the $\ZZ^m$-multigraded generic initial ideal of $J$ is generated by $B_{m,n}$ implies that the $K$-vector space $\ini_< (\phi(J_m))$ is generated by $B_{m,n}$. Hence the initial algebra of the $K$-algebra $\phi(K[T])$ contains $K[B_{m,n}]$. It follows that  
 $$H(K[T],u)= H(\phi(K[T]),u)\geq H(K[B_{m,n}],u)$$
  for all $u\in \ZZ^m$ and, in particular, $u\in \ZZ$. Hence
  $$\dim K[B_{m,n}] \leq \dim K[T] \leq \dim G(m,n)=m(n-m)+1.$$
  We finally observe that the subset 
  $$B_{m,n}^0=\bigcup_{i=1}^m  \left \{  \left(\prod_{k=1}^m X_{k1} \right) X_{ij}/X_{i1} \  : j=1,2,\dots, n-m+1 \right\}$$  
 of $B_{m,n}$ consists of algebraically independent elements and has $m(n-m)+1$ elements. 
 Hence $\dim K[B_{m,n}]\geq m(n-m)+1$, concluding the proof. 
\end{proof}

\section{Are the minors a \Sagbi{} basis?}
\label{AreThey?}

Denote by $ \cM_t$ the set of $t$-minors of  $X_{m\times n}$. As we have already recalled, $ \cM_m$  is  a  \Sagbi{}  basis of $G(m,n)$ with respect to the diagonal order. Furthermore, for $1<t<m$ the set $\cM_t$ is not a  \Sagbi{} basis of  $A_t(m,n)$ with respect to the diagonal order. The most immediate argument for this claim is that $\dim A_t = mn$ if $1 <t <m$  \cite[Prop. 6.4.1]{BCRV}, but not all variables in the matrix $X$ occur in diagonals. In the following the characteristic of the field is irrelevant.

\begin{theorem}\label{AtSagbi}	
There exists a monomial order on $K[X_{m\times n}]$ for which $\cM_t$ is a \Sagbi{} basis of $A_t(m,n)$ if and only if  
\begin{enumerate}
\item $m=t$, i.e., $A_t(m,n) = G(m,n)$ and  $n \ge t$;
\item $m = n = t+1$;
\item $t=1$.
\end{enumerate}
\end{theorem}

Case (3) is trivial,  and we have just recalled (1). As a next step we prove the existence of a suitable order in case (2).

\begin{proposition}\label{submax}
The set $ \cM_{m-1}$ is  a \Sagbi{} basis of $A_{m-1}(m,m)$ with respect to the Lex  order associated to the total order 
$$
X_{11} > X_{22} > \cdots > X_{mm} > \cdots \mbox{ the remaining variables in some order}.
$$
\end{proposition} 

\begin{proof}  The elements of $ \cM_{m-1}$ in $K[X_{m\times m}]$ are algebraically independent.  So  $ \cM_{m-1}$ is  a \Sagbi{} basis with respect to a monomial order if and only if the  initial monomials of the elements in $ \cM_{m-1}$ are algebraically independent. Therefore we must check  this for the given Lex order. Given 
 $i,j\in [m]$,  for the $(m-1)$-minor  $\mu_{ij}$ that does not use row $i$ and column $j$ we have 
\begin{equation}
\label{inform} 
\ini(\mu_{ij}) =X_{ji}  \Delta / X_{ii}X_{jj}
\end{equation} 
where  $\Delta=\prod_{k=1}^m  X_{kk}$. Consider the field extension $L$ of $K$ generated by all the initial monomials $\ini(\mu_{ij})$. We must show that it has transcendence degree $m^2$ over $K$.  We have
$$
\Delta^{m-1} = \prod_{i=1}^m \ini(\mu_{ii}).
$$
The equation shows that $\Delta$ is algebraic over $L$. Therefore we can extend $L$ by $\Delta$ to $L'$ without changing the transcendence degree. We claim that $L'$ contains all indeterminates $X_{ij}$, $i,j\in[m]$ and therefore is the full fraction field of $K[X]$. 
Since
$$
X_{ii} = \Delta/\ini(\mu_{ii})
$$
it follows that the indeterminates $X_{ii}$, $i\in [m]$ are in $L'$. But then Eq.(\ref{inform}) 
 implies that $X_{ji}$, $j,i\in [m]$, $j\neq i$, also belong to $L'$, and we are done.
 \end{proof}

  \begin{definition} If $m\neq n$ we denote  by  $G_{m\times n}$ be the product $S_m\times S_n$ of the symmetric groups  $S_m$ and $S_n$. Furthermore we denote by $G_{m\times m}$ the  semi direct product of $S_m\times S_m$ and the cyclic group of two elements. 
  \end{definition} 
  
  The group $G_{m\times n}$ acts on the polynomial ring $K[X_{m\times n}]$ with $S_m$ and $S_n$ permuting the rows and columns of $X$ and, when $m=n$,  the cyclic group of order $2$ transposing the matrix. 

  A crucial point in the proof that $\cM_t$ is never a \Sagbi{} basis in all cases different from (1), (2) and (3) in Theorem \ref{AtSagbi} is that the selection of initial monomials in the proof of Proposition \ref{submax} is  unique up to the action of $G_{m\times m}$.

 The group $G_{m\times n}$ acts on the set $\cM_t$ (up to sign) and hence also on the Newton polytope of  $\cM_t$.
 
 As done in \cite{SZ}, when dealing with polynomials in $K[X_{m\times n}]$ it is natural to represent the exponent vector of a monomial  as a $m\times n$ \emph{exponent matrix} with entries in $\NN$. Similarly a weight vector on $K[X_{m\times n}]$ can be given as a $m\times n$ matrix with entries in $\NN$. 
 
 \begin{remark}\label{magicsquare} 
 Every matching of $\cM_t$ and, in particular, every vertex  of the Newton polytope $\New(\cM_t)$ of $\cM_t$,  when  represented as a $m\times n$ matrix with entries in $\NN$ has the following properties: 
 \begin{itemize} 
\item[(1)]  each row sum equals  
$${m-1 \choose t-1}{n \choose t},$$
\item[(2)]  each column sum  equals
$${m \choose t}{n-1 \choose t-1}.$$
\end{itemize} 
This is easily seen because the minors are multi homogeneous with respect to the $\ZZ^m\times \ZZ^n$ grading given by $\deg X_{ij}=(e_i, e_j)\in \ZZ^m\times \ZZ^n$. In other words,  the number in (1)  is the number of $t$-minors involving a given row.  Similarly for columns. In the square case $m = n$ it is a  ``magic square" in the sense of Stanley \cite{Stan}.     
\end{remark} 

 Inspection of the proof of Proposition \ref{submax} shows that the coherent matching constructed there has the exponent matrix 
 $$
 Q_m= d_mI_m +E_m 
 $$ 
 where 
 
 \begin{itemize}
 \item[(i)] $I_m$ is the identity matrix of size $m\times m$, 
 \item[(ii)] $d_m=(m-1)^2-1$  and 
 \item[(iii)] $E_m$ is  a matrix of size  $m\times m$ with all entries $1$.  
 \end{itemize}
 
 As an example for $m=3$ one has: 
$$
Q_3 =
 3\begin{pmatrix}
	1 &0 &0 \\
	0 &1 &0 \\
	0 &0 &1
\end{pmatrix}
+
\begin{pmatrix}
	1 &1 &1 \\
	1 &1 &1 \\
	1 &1 &1
\end{pmatrix}
=\begin{pmatrix}
	4 &1 &1 \\
	1 &4 &1 \\
	1 &1 &4
\end{pmatrix}
.$$
 
 A matrix in  $\NN^{m\times n}$ is said to have \emph{full support} if it has only positive entries. 

\begin{theorem}\label{unique} 
For $m\geq 2$ consider the set $\cM_{m-1}$ of $(m-1)$-minors of  $X_{m\times m}$  and the Newton polytope $\New(\cM_{m-1})$, i.e the Newton polytope of the product of the elements in $\cM_{m-1}$.  Up to  the $G_{m\times m}$ action, $Q_m$ is the only vertex of $\New(\cM_{m-1})$ with full support.  In other words, up to  the $G_{m\times m}$ action, the only coherent matching of $\cM_{m-1}$  involving all variables is the one given in Proposition  \ref{submax}.
\end{theorem} 

\begin{proof}
The statement is trivially true for $m=2$ so that we can assume $m\ge 3$ in the following.
Let $P$  be the exponent matrix of an arbitrary matching of  $\cM_{m-1}$ and assume that $P$ has full support. It is enough to prove: 
\smallskip 

\noindent  {\bf Claim.} \enspace
The matrix $P$  is in the polytope spanned by $Q_m$ and its conjugates under the action of $G_{m\times m}$. 
\smallskip 
  
As observed in Remark \ref{magicsquare} $P$ is a ``magic square'': all row and column sums have the constant value $m(m-1)$ and, by assumption,  all entries of $P$ are  positive integers. 

 To check the claim, and since $E_m$ is fixed by $G_{m\times m}$,  we may translate both $P$ and $Q_m$ by subtracting $E_m$. In this way  $P-E_m$ is a matrix with non-negative entries with row and column sums equal  to $d_m$ and $Q_m-E_m = d_mI_m$. 
 We may as well multiply both $P-E_m$ and $d_mI_m$ with $d_m^{-1}$. It follows that  $d_m^{-1}(P-E_m)$ is a matrix with non-negative entries  and row and column sums equal to $1$, hence a point of the Birkhoff polytope $B_m$  (see \cite{Stan}) with $m$ rows.   The vertices of $B_m$ are know to be the permutation matrices by the Birkhoff--von Neumann theorem \cite[p. 32]{Stan} , i.e.  the orbit of $I_m$  under the action of $G_{m\times m}$.  This concludes the proof.  \end{proof}

As an immediate consequence of  Proposition \ref{submax} and Theorem \ref{unique} we have: 

\begin{corollary} 
Up to the $G_{m\times m}$ action, the only coherent matching that makes 
$\cM_{m-1}$ a \Sagbi{} basis of $A_{m-1}(m,m)$ is the one of Proposition \ref{submax}. 
\end{corollary}   

\begin{proof} The elements of $\cM_{m-1}$ in $K[X_{m\times m}]$ are algebraically independent. Hence  a coherent matching that makes $\cM_{m-1}$ a \Sagbi{}  basis must involve all the variables. In other words, the corresponding vertex of $\New(\cM_{m-1})$ must have full support. But then by Theorem \ref{unique} we know that such a coherent matching must equal to that of Proposition \ref{submax} up to the $G_{m\times m}$ action. 
\end{proof}

\begin{remark}
We illustrate Theorem \ref{unique} with some data  on the polytope $\New(\cM_{m-1})$ with   $m=3$ and $m=4$. 

(a) Consider $ \cM_2$ as a subset of $K[X_{3\times 3}]$. The vertices of  the Newton polytope $\New(\cM_2)$   as well as the their orbits under $G_{3\times 3}$ can be computed by Normaliz \cite{Nmz}.  In total there are $102$ vertices and  $5$ orbits, represented by  the exponent  matrices in $\ZZ^{3\times 3}$ given in the  Table \ref{3x3}. 
\begin{table}
 $$
 \begin{array}{ccccc} 
         (1) & (2) & (3) &  (4) &  (5)   \\ 
         \hline 
         \\
   \begin{pmatrix}
        4 &2 &0 \\
        2 &2 &2 \\
        0 &2 &4
    \end{pmatrix} & 
  \begin{pmatrix}
        4 &2 &0 \\
        2 &1 &3 \\
        0 &3 &3
    \end{pmatrix} & 
     \begin{pmatrix}
        3 &2 &1 \\
        2 &0 &4 \\
        1 &4 &1
    \end{pmatrix} &
  \begin{pmatrix}
        4 &1 &1 \\
        1 &4 &1 \\
        1 &1 &4
    \end{pmatrix} &
  \begin{pmatrix}
        0 &3 &3 \\
        3 &0 &3 \\
        3 &3 &0
         \end{pmatrix}
         \end{array}
 $$
 \vspace*{1ex}
  \caption{Representatives of orbits of vertices of $\New(\cM_2)$ in $3\times 3$}\label{3x3}
\end{table}
 They correspond bijectively to the coherent matchings of $\cM_2$ up to $G_{3\times 3}$. As predicted by Theorem \ref{unique}   there is only one vertex with full support,  (4)  in Table \ref{3x3},  and it is exactly~$Q_3$. 
 
 (b) For $m=4$ the Newton polytope $\New(\cM_{m-1})$ has $77328$ vertices. The number was predicted by evaluating random linear forms on the Newton polytope and confirmed by Normaliz within $\sim 4$ days. The vertices decompose into $98$ orbits under the action of $G_{4\times 4}$ and only the orbit of $Q_4$ has full support. 
\end{remark} 
  
 A simple tool is the following evident lemma that will allow us to pass to submatrices.

\begin{lemma}\label{passage}
Let $R$ be a polynomial ring over a field endowed with a monomial order. Let $A$ and $S$ be a subalgebras of $R$  with $S$  generated by a subset of the variables. Let   $\cF \subset A$.  Suppose that  for all $F\in \cF$   with $\ini(F)\in S$ one has $F\in S$. 
\begin{enumerate}
\item If $\cF$ is a \Sagbi{} basis of $A$, then $\cF\cap S$ is a \Sagbi{} basis of $A\cap S$.
\item Assume  $A$ is graded and the elements in $\cF$ are homogeneous.  If $K[\ini(\cF)]_i = \ini(A)_i$ for every $i=1,\dots,k$ then $K[\ini(\cF\cap S)]_i = \ini(A\cap S)_i$ for every $i=1,\dots,k$.
\end{enumerate}
\end{lemma}
  
We are ready to prove  Theorem \ref{AtSagbi}.

\begin{proof}[Proof of Theorem \ref{AtSagbi}] 
It remains to prove that in all cases different from (1), (2) and (3)  there is no monomial order such that $\cM_t$ is a \Sagbi{} basis. 
Suppose, by contradiction,  that there exists a monomial order $<$ on $K[X_{m\times n}]$ such that $\cM_t$ is a \Sagbi{} basis. Being a  \Sagbi{} basis is compatible with the restriction to a submatrix, as follows from Lemma \ref{passage}: if the initial monomial of a minor ``lives'' in the submatrix, then the minor lives in the submatrix. So we may assume right away that $t = m-1$ and the format of our matrix $X$ is $m \times (m+1)$. In the following we will denote by $X^{(j)}$ the $m\times m$ matrix obtained from $X$ by removing column $j$.

By Theorem \ref{unique}  for every $j=1,\dots, m+1$ the exponent  matrix $M^{(j)}$ of the product of the monomials in the  matching restricted to the submatrix $X^{(j)}$ must be a conjugate of $Q_m$. This will lead to a contradiction.  Indeed if an  entry of $M^{(j)}$  has a value $>1$, this value must be $(m-1)^2$. Moreover, if we remove a further column $k\neq j$ from $X^{(j)}$, then the remaining exponent  matrix has entries $> 1$ at all indeterminates that had degree $(m-1)^2$ in $M^{(j)}$  and do not belong to column $k$: we say that  $X^{(k)}$ ``inherits'' all variables $X_{uv}$ from  $X^{(j)}$ whose degree is $(m-1)^2$ and have $v\neq k$.

We start from  $X^{(m+1)}$, and may assume that $X_{11},\dots, X_{mm}$ have degree $(m-1)^2$ with respect to it. 

Then $X^{(m)}$ inherits $X_{11},\dots, X_{m-1,m-1}$ from $X^{(m+1)}$, and therefore $X_{m,m+1}$ must have degree $(m-1)^2$  in $M^{(m)}$  as well.

Take  $X^{(m-1)}$. From $X^{(m+1)}$ it inherits $X_{11},\dots, X_{m-2,m-2}$ and $X_{mm}$. So $X_{m-1,m+1}$ must have degree $(m-1)^2$  in  $M^{(m-1)}$ as well. But from $X^{(m)}$ we also get that $X_{m,m+1}$  has degree $(m-1)^2$ in $M^{(m-1)}$, and this is a contradiction: exactly $m$ variables must have degree $(m-1)^2$ in $M^{(m-1)}$.
\end{proof}

\begin{remark}  In the proof of Theorem \ref{AtSagbi} for $m=3$, after the reduction to the $3\times 4$ case,  one can proceed also  by ``brute force" using Hilbert series.  Consider the product $F$ of the elements in $\cM_2$ and compute the vertices of its Newton polytope and their orbits by Normaliz\cite{Nmz}. There are $3624$ vertices and  $29$ orbits  under $G_{3\times 4}$.  Among  the $29$ orbits  only the $5$ in Table \ref{3x4} have full support. 
\begin{table}
\begin{small}
$$
 \begin{array}{ccccc} 
         (1) & (2) & (3) &  (4) &  (5)   \\ 
        \hline 
         \\ 
v_1= \begin{pmatrix}
2 & 2 & 2 & 6 \\
6 & 1 & 3 & 2 \\
1 & 6 & 4 & 1 \\
\end{pmatrix}
    &  \! \!  \! \! 
 \begin{pmatrix}
3 & 6 & 2 & 1 \\
1 & 2 & 6 & 3 \\
5 & 1 & 1 & 5 \\
\end{pmatrix}
    &  \! \!  \! \! 
 \begin{pmatrix}
6 & 3 & 1 & 2 \\
2 & 5 & 4 & 1 \\
1 & 1 & 4 & 6 \\
\end{pmatrix}
     & \! \!  \! \! 
 \begin{pmatrix}
3 & 2 & 1 & 6 \\
3 & 6 & 2 & 1 \\
3 & 1 & 6 & 2 \\
\end{pmatrix}
    & \! \!  \! \! 
 \begin{pmatrix}
6 & 4 & 1 & 1 \\
1 & 1 & 4 & 6 \\
2 & 4 & 4 & 2 \\
\end{pmatrix}
         \\ \\
w_1=\begin{pmatrix}
1 & 0 & 2 & 3 \\
3 & 0 & 3 & 2 \\
0 & 2 & 2 & 0 \\
\end{pmatrix}
       &  \! \!  \! \! 
 \begin{pmatrix}
2 & 3 & 0 & 1 \\
1 & 1 & 3 & 3 \\
3 & 0 & 0 & 3 \\
\end{pmatrix}
   &  \! \!  \! \! 
 \begin{pmatrix}
2 & 1 & 1 & 0 \\
1 & 2 & 3 & 0 \\
0 & 0 & 3 & 3 \\
\end{pmatrix}
    & \! \!  \! \! 
\begin{pmatrix}
0 & 1 & 0 & 3 \\
0 & 3 & 1 & 0 \\
0 & 0 & 3 & 1 \\
\end{pmatrix}
  &  \! \!  \! \! 
\begin{pmatrix}
5 & 3 & 2 & 1 \\
1 & 0 & 4 & 5 \\
2 & 2 & 4 & 2 \\
\end{pmatrix}  
\\ \\
h_1=(1,6,11,5) 
& 
(1,6,10,3)  
& 
(1,6,10,4)
&
(1,6,12,7)   
&
(1,6,9,4)   
         \end{array}
$$
\end{small}
 \vspace*{1ex}
\caption{Vertices with full support of $\New(\cM_2)$ in  $3\times 4$}\label{3x4}
\end{table}
There the first row shows the five  full support vertices, while the second row shows a weight corresponding to the vertex. 

For example, the coherent matching of $\cM_2$ associated to vertex (1) is given by the  initial monomials of the $2$-minors with respect to the weight $w_1$. They are 

\begin{equation*}
\label{fullSup} 
\begin{array}{cccccc}
X_{12}X_{21},  &  X_{13}X_{21},  &  X_{14}X_{21}, & X_{12}X_{23}, &  X_{14}X_{22}, &  X_{14}X_{23}, \\
X_{11}X_{32},  & X_{11}X_{33},  & X_{14}X_{31},  &X_{13}X_{32},  & X_{14}X_{32}, & X_{14}X_{33}, \\
X_{21}X_{32}, & X_{21}X_{33}, & X_{21}X_{34}, &  X_{23}X_{32}, & X_{24}X_{32}, & X_{24}X_{33}.
\end{array}
\end{equation*}
As said already,  the weight is only a tool to get these initial monomials: they are uniquely determined by the vertex of the Newton polytope.

Now the Hilbert series of the $K$-algebra generated by these $18$ monomials is easily computed as 
$$
(1 + 6z + 11z^2 + 5z^3) / (1-z)^{12}
$$
and $h$-vector $h_1=(1,6,11,5)$, 
which is different from the Hilbert series 
$$
(1 + 6z+ 15z^2 + 10z^3) / (1-z)^{12}
$$
of $A_2(3,4)$. This shows that the $2$-minors are not a \Sagbi{} basis with respect to any order compatible with the given coherent matching. 
The same computation can be repeated for the other four vertices. The corresponding  $h$-vectors are in the third row of in Table \ref{3x4}. 
This completes the alternative proof that no order can make the $2$-minors of $X_{3\times 4}$ a \Sagbi{} basis. 
\end{remark} 

\section{Universal \Sagbi{} bases for $A_2(3,3)$,  $G(3,6)$ and $G(3,7)$} 
\label{univSag}

Whereas the explicit computations in the previous section only served heuristic or illustrative purposes, they are crucial for the results in this section. All proofs are by computation. 

A universal  \Sagbi{} basis for an algebra $A$ is a set of elements that are a  \Sagbi{} basis for all orders. 
In this section we discuss universal  \Sagbi{} bases for two algebras $A_2(3,3)$, $G(3,6)$ and $G(3,7)$. We start with the first: 
 
  \begin{theorem} 
  \label{univSagbi3x3} 
  Let $\Delta$ be the determinant of $X_{3\times 3}$.  Then the set 
  $$\cU=\cM_2\cup \{ X_{ij}\Delta : i,j\in [3]\}$$ is a universal \Sagbi{}  basis of  $A_2(3,3)$. 
  More precisely, for every order a \Sagbi{} basis of  $A_2(3,3)$ is obtained by adding to $\cM_2$ at most three  elements of type  $X_{ij}\Delta$. 
  \end{theorem} 
  
  \begin{proof} The fact that $\cU$ is a subset of $A_2(3,3)$ is  a special case of \cite[Lemma 3.5.6]{BCRV}. 
  Consider  the coherent matchings  of $\cM_2\cup\{\Delta\}$.  With   Normaliz\cite{Nmz} we have checked   that the Newton polytope of the product of the polynomials  $\cM_2\cup\{\Delta\}$  has  $108$ vertices in total and  $5$ distinct orbits  (see Table \ref{3x3Delta})  up to the action of $G_{3\times 3}$.
  \begin{table}
  	$$
  	\begin{array}{ccccc} 
  		(1) & (2) & (3) &  (4) &  (5)   \\ 
  		\hline 
  		\\
  		\begin{pmatrix}
  			5 &2 &0 \\
  			2 &3 &2 \\
  			0 &2 &5
  		\end{pmatrix} & 
  		\begin{pmatrix}
  			5 &2 &0 \\
  			2 &1 &4 \\
  			0 &4 &3
  		\end{pmatrix} & 
  		\begin{pmatrix}
  			4 &2 &1 \\
  			2 &0 &5 \\
  			1 &5 &1
  		\end{pmatrix} &
  		\begin{pmatrix}
  			5 &1 &1 \\
  			1 &5 &1 \\
  			1 &1 &5
  		\end{pmatrix} &
  		\begin{pmatrix}
  			0 &4 &3 \\
  			3 &0 &4 \\
  			4 &3 &0
  		\end{pmatrix}
  	\end{array}
  	$$
  	\vspace*{1ex}
  	\caption{Representatives of orbits of vertices of $\New(\cM_2\cup{\Delta})$ in $3\times 3$}\label{3x3Delta}
  \end{table}   
 Each vertex in Table \ref{3x3Delta} is obtained as the sum of the  vertex of $\New(\cM_2)$ in Table \ref{3x3}  (indexed  with the same number) and a monomial of $\Delta$.     
  
 Case (4) is the easiest: it has been proved in Lemma \ref{submax} that $\cM_2$ is a \Sagbi{} basis for it.  Next we analyze case (1) in Table \ref{3x3Delta} (defined by a diagonal order). It has two entries $0$ in the support at positions $(1,3)$ and $(3,1)$. The  associated coherent matching of $\cM_2\cup\{\ \Delta\}$  is 
  $$
T=X_{11}X_{22},\  X_{11}X_{23}, \  X_{12}X_{23}, \ X_{11}X_{32}, \   X_{11}X_{33}, \   X_{12}X_{33}, \   X_{21}X_{32}, \   X_{21}X_{33}
$$ 
together with  the initial $D=X_{11}X_{22}X_{33}$ of $\Delta$.
The variables $X_{13}$ and $X_{31}$ do not appear in $T$. But $X_{13}D$ and $X_{31}D$ are both initial monomials of elements of $\cU$. So they are algebraically independent over $K[T]$ because they involve variables $X_{13}$ and $X_{31}$ that are not present in $T$.  Now, to prove that $\cU$ is a \Sagbi{} basis of $A_2(3,3)$ with respect to any monomial order compatible with the coherent matching  (1) it is enough to prove that  
$$
C=K[T, X_{13}D, X_{31}D]
$$ 
has the Hilbert series of $A_2(3,3)$, that is,  $1/(1-z)^9$. The Hilbert series of $C$ is that of $K[T]$ divided by $(1-z^2)^2$. Summing up, it is enough to show that the Hilbert series of $K[T]$ is 
$$
\frac{(1-z^2)^2}{(1-z)^9}
$$
and this is easily checkable by direct computation: the defining ideal of $K[T]$ is a complete intersection of $2$ quadrics. 

The same argument can be applied to the remaining $3$ cases. It is important  that the initial monomial of $\Delta$ is a product of indeterminates that already appear in the corresponding exponent matrix in Table \ref{3x3} with a positive degree. This allows us to use the same argument for algebraic independence as in case (1). In each case the elements that we have to add to $\cM_2$ to get a  \Sagbi{} are of the form $X_{i,j}\Delta$ where the $(i,j)$-entry of vertex in Table \ref{3x3Delta} equals $0$. So we have to add at most three at each time.  
\end{proof} 

\begin{remark}\label{not_unique}
In the  proof of  Theorem  \ref{univSagbi3x3}  we have observed that the vertices of the Newton polytopes of $\cM_2$ and of $\cM_2\cup \{\Delta\}$ have the same number of orbits under $G_{3\times 3}$. This seems to suggest that the initial monomial of $\Delta$ is uniquely determined by the initial monomials of the $2$-minors. This is true for the cases (1)--(4), but wrong in case (5).  Indeed given the coherent matching on $\cM_2$ that corresponds to (5) in Table \ref{3x3} there are two potential initial monomials of $\Delta$: either $X_{12}X_{23}X_{31}$ or   $X_{21}X_{32}X_{13}$ and both are possible.  
With the first choice we get vertex (5) in Table \ref{3x3Delta} while with the second choice we get 
$$
 \begin{pmatrix}
	0 &3 &4 \\
	4 &0 &3 \\
	3 &4 &0
\end{pmatrix}
$$
which is however  in the $G_{3\times 3}$ orbit  of the vertex  (5) of Table \ref{3x3Delta}.  
Despite of the same number of orbits, $\New(\cM_2\cup{\Delta})$ has $108$  vertices in total, whereas $\New(\cM_2)$ has only $102$. The orbits of case (5) of Table \ref{3x3} has $6$ elements while that of (5) of Table \ref{3x3Delta}  has $12$ elements.
\end{remark}

\begin{remark}\label{product}
(a) For a diagonal order $A_t(m,n)$ has a well-understood \Sagbi{} basis of products of minors if $\chara K =0$ or $\chara K > \min(t, m-t, n-t)$.  In fact, the initial algebra is precisely described in \cite[Lemma 6.4.2]{BCRV}, and by a theorem of Varbaro \cite[Thm. 6.4.10]{BCRV}, the maximum (relative) degree in a minimal system of generators of $\ini(A_t(m,n))$ is $\leq m-1$. The assumption on characteristic is always satisfied for $m= n$, $t = m-1$.   
In the case $m = n = 4$, $t = 3$,  the initial algebra has two generators in degree $3$. Indeed generators of the initial algebra  are given by initials of product of minors of shape
\begin{enumerate}
	\item[(i)]  (3), i.e the $3$ minors,
	\item[(ii)]  (4,2), i.e. the product  of the determinant $\Delta$ and a  $2$-minors, and
	\item[(iii)]  (4,4,1) i.e. polynomials of type $X_{ij} \Delta^2$.
\end{enumerate}
The explicit computation shows  that $10$ elements of type (ii) and only $X_{14}\Delta^2$ and $X_{41}\Delta^2$ of type (iii) are actually needed in the \Sagbi{} basis. In particular  every universal  \Sagbi{} basis of  $A_3(4,4)$ must contain  elements of degree $3$. 

(b) It turns out that already in  $A_3(4,4)$ the products of minors that belong to the algebra are not a universal   \Sagbi{} basis. Indeed the product of minors that are in $A_3(4,4)$  are the one  listed above or further multiples.  With respect to the Lex order associated to the order of the variable 
$$
\begin{array}{l}
	X_{11}> X_{13}> X_{12}> X_{22}> X_{14}> X_{23}> X_{21}> X_{24}>  \\ 
	X_{31}> X_{32}> X_{33}> X_{34}> X_{41}> X_{42}> X_{43}> X_{44}
\end{array} 
$$
they do not form a \Sagbi{} basis. Already in degree $3$ one needs extra elements. 
\end{remark}

Now consider the algebra $G(3,6)$. The products of minors in the algebra are not a universal \Sagbi{} basis: the only products of minors in $G(3,6)$ are the products of maximal minors and the maximal minors do not form a universal \Sagbi{} basis. The way out is to take the polynomial 
$$F=[1,2,3][4,5,6]-[1,2,4][3,5,6]$$
and its orbit 
$$
O_F=\bigl\{ [a,b,c][d,e,f]-[a,b,d][c,e,f]  : \{a,b,c,d,e,f\}=[6] \bigr\} \mod \pm 1
$$
under $S_6$ permuting the columns. The polynomial $F$ has $48$ terms and, up to sign, $F$ is fixed by $48$ permutations of columns. Hence, up to sign, 
 $O_F$ consists of $15$ different polynomials.   Our goal is to prove:

  \begin{theorem} 
  \label{univSagbi3x6} 
  The set  $\cM_3\cup O_F$ is a universal \Sagbi{}  basis of  $G(3,6)$. 
  More precisely, for every order a \Sagbi{} basis of  $G(3,6)$ is obtained by adding at most one  element of $O_F$ to $\cM_3$.   
  \end{theorem} 
  
  \begin{proof} 
 For any coherent matching $T$ of  $\cM_3$ we will check that either
 \begin{itemize} 
 \item[(i)] $T$ gives already a toric algebra with the  Hilbert series 
 \begin{equation}
\label{HS36}
(1 + 10z + 20z^2 + 10z^3 + z^4) / (1-z)^{10}
\end{equation} 
of $G(3,6)$, or 
\item[(ii)] there exists  $G\in O_F$ such that all the extensions of $T$ to a coherent matching of $\cM_3\cup \{G\}$ gives  a toric algebra with the  Hilbert series of $G(3,6)$.
\end{itemize} 
 
To do this, we analyze the Newton polytope $\New(\cM_3)$  associated to  the product of the $3$-minors.  More precisely we will analyze the faces  of $\New(\cM_3)$ that correspond to the standard form of the four types  identified by  \cite{SZ}. Since the set $\cM_3\cup O_F$ is fixed (up to sign) by the group $G_{3\times 6}$ this is enough to treat all the coherent matchings of $\cM_3$.      We start with: 
  
 {\bf  (Type 1)} These are the vertices of $\New(\cM_3)$ that, up to symmetry, have the form:  
  
 $$ \begin{pmatrix}
  			* &0 &0  & * &* & * \\
  			0 &* &0 & * &* & * \\
  			0 &0 &* & * &* & * \\
  		\end{pmatrix} 
.$$		
They correspond to coherent matchings that do not involve the variables 

$$\FV=\{X_{12}, X_{13}, X_{21}, X_{23}, X_{31}, X_{32}\} $$ associated to the entries $0$. We call the elements in $\FV$ the ``forgotten variables". In practice, we may replace in the above matrix the $*$ with variables,
\begin{equation}
\label{Type1mat}
 X^{(1)}=\begin{pmatrix}
  			X_{11} &0 &0  & X_{14} &X_{15} & X_{16} \\
  			0 &X_{22} &0 & X_{24} &X_{25} & X_{26}\\
  			0 &0 &X_{33} & X_{34} &X_{35} & X_{36} \\
  		\end{pmatrix} 
\end{equation} 	
 then take the product of the set $\cM_3^{(1)}$ of the $3$-minors of $X^{(1)}$ and compute the vertices of the Newton polytope up to the subgroup (of order $36$) of $G_{3\times 6}$  of rows and columns permutations fixing the shape above.  It turns out that  there are  $108$ vertices and $5$ orbits. Not surprisingly  we get the same numbers that we have already seen the proof  Theorem \ref{univSagbi3x3}, see \cite[Prop.3.11]{SZ} for an explanation.  Orbit representatives are the following: 
 $$
  	\begin{array}{ccc} 
  		(1) & (2) & (3)   \\ 
  		\hline 
  		\\
  		\begin{pmatrix}
  10& 0& 0 &1 & 3 & 6 \\
  0  &10 & 0 & 6& 3 & 1 \\
  0 & 0 & 10 &3& 4 & 3
  		\end{pmatrix} & 
  		\begin{pmatrix}
  10& 0 & 0 & 2 & 6 &  2\\
 0 & 10 & 0 & 3 & 1 & 6 \\
  0 & 0  & 10 & 5 & 3& 2 
  		\end{pmatrix} & 
 \begin{pmatrix}
   10& 0& 0& 4& 5& 1\\
  0& 10& 0& 5& 1& 4 \\
  0& 0& 10& 1& 4& 5
  		\end{pmatrix}  		
  	\end{array}
  	$$

$$
  	\begin{array}{ccc} 
  		(4) & (5)    \\ 
  		\hline 
  		\\
  		\begin{pmatrix}
   10& 0& 0& 5& 3& 2 \\
  0& 10& 0& 4& 1& 5 \\
  0& 0& 10& 1& 6& 3 
  		\end{pmatrix} & 		
		\begin{pmatrix}
  10& 0& 0& 6& 2&2 \\ 
  0& 10& 0& 2& 6& 2 \\
  0& 0&10& 2& 2& 6 
  		\end{pmatrix}  
  	\end{array}
  	$$
The vertices (1)-(4) correspond to coherent matchings whose associated toric algebra has Hilbert series (\ref{HS36})  while (5) gives the matching 
\begin{equation} 
\label{mathing5}
\begin{array}{cccccc}
T=& X_{11}X_{22}X_{33}, & X_{11}X_{22}X_{34}, & X_{11}X_{22}X_{35}, &X_{11}X_{22}X_{36}, & X_{11}X_{24}X_{33} \\ 
& X_{11}X_{25}X_{33}, & X_{11}X_{26}X_{33}, & X_{11}X_{25}X_{34}, & X_{11}X_{24}X_{36}, &X_{11}X_{25}X_{36}, \\
& X_{14}X_{22}X_{33}, & X_{15}X_{22}X_{33}, & X_{16}X_{22}X_{33}, & X_{14}X_{22}X_{35}, & X_{14}X_{22}X_{36},\\ 
& X_{15}X_{22}X_{36}, & X_{14}X_{25}X_{33}, & X_{14}X_{26}X_{33}, & X_{16}X_{25}X_{33}, & X_{14}X_{25}X_{36}
\end{array}
\end{equation} 
and Hilbert series 
\begin{equation}
\label{NoHS36}
(1 + 10z + 19z^2 + 8z^3) / (1-z)^{10}.
\end{equation} 
Hence for coherent matchings  associated to (1)-(4)  we are in case (i).   On the other hand (5) does not give the correct Hilbert series.  
As predicted by Theorem \ref{dimok}, we see in  (\ref{NoHS36}) that  for (5) we have already reached the correct Krull dimension.  Hence  when extending the coherent matching $T$  to a coherent matching of $\cM_3\cup\{G\}$, for any  $G$ in $G(3,6)$ none of the forgotten variables can appear.  
In other words,  the forgotten variables $\FV$   for Type 1 are  forgotten forever. 
We take the following element of $O_F$:
$$
G=[1, 4, 2][5, 3,6]-[1, 4, 5][2, 3,6].
$$
Expanding the minors  we have: 
\begin{equation}
\label{thisisG} 
G=G_0 \mod \FV
\end{equation} 
with 
$$G_0=X_{11}X_{22}X_{33}X_{15}X_{26}X_{34} - X_{11}X_{22}X_{33}X_{16}X_{24}X_{35}.$$ 

Hence every extension of $T$ to a coherent matching of $\cM_3\cup \{G\}$ can only involve one the two monomials of $G_0$. 
Now the last step:  we consider   the coherent matchings of  $\cM_3^{(1)}\cup \{G_0\}$ that extend the coherent matching (\ref{mathing5}).  This amounts to taking the Newton polytope of  the product of the elements in $\cM_3^{(1)}\cup \{G_0\}$  and to identify the vertices that are compatible with (5). It turns out that, up to symmetry, there is only one such vertex. It  corresponds to the coherent matching 
$$T_1=T, X_{11}X_{22}X_{33}X_{15}X_{26}X_{34}.$$ 
Finally we compute the Hilbert series of the toric algebra associated to $T_1$ and verify  that it has Hilbert series (\ref{HS36}).  So we have verified (ii). We finally note that (5) is the only orbit representative of a vertex of Type 1 with only even coordinates.  

The same scheme applies to the other $3$ types. We just record below the salient data. 
 
 {\bf  (Type 2)}. They correspond to coherent matchings that do not involve the variables 
 $$\FV=\{X_{21}, X_{31}, X_{32}, X_{15}, X_{16}, X_{26}\} $$
 and hence vertices of the Newton polytope  
  $$ \begin{pmatrix}
  			* & * &*  & * & 0 &  0 \\
  			0 &* &* & * &* &  0  \\
  			0 &0 &* & * &* & * \\
  		\end{pmatrix} 
.$$		
 Here 
\begin{equation}
\label{Type2mat}
 X^{(2)}=\begin{pmatrix}
  			X_{11} &X_{12} &X_{13}   & X_{14} & 0  &  0  \\
  			0 &X_{22} &X_{23}  & X_{24} &X_{25} & 0\\
  			0 &0 &X_{33} & X_{34} &X_{35} & X_{36} \\
  		\end{pmatrix} 
\end{equation} 	
and the $\New(\cM_3^{(2)})$ has $80$ vertices and $22$ orbits under the subgroup (of order $4$) of $G_{3\times 6}$ fixing the shape above. Only one orbit does not correspond to a \Sagbi{} basis and it has orbit representative
$$
\begin{pmatrix}
  10& 2& 6& 2& 0&0 \\ 
  0& 8& 2& 2& 8& 0 \\
  0& 0&2& 6& 2& 10 
  \end{pmatrix}  
  $$ associated coherent matching 
  \begin{equation} 
\begin{array}{cccccc}
T=& X_{11}X_{22}X_{33}, &X_{11}X_{22}X_{34}, &X_{11}X_{22}X_{35}, & X_{11}X_{22}X_{36}, &X_{11}X_{23}X_{34}, \\
&X_{11}X_{25}X_{33}, & X_{11}X_{23}X_{36}, & X_{11}X_{25}X_{34}, & X_{11}X_{24}X_{36}, & X_{11}X_{25}X_{36}, \\
& X_{13}X_{22}X_{34}, &X_{13}X_{22}X_{35},& X_{13}X_{22}X_{36},& X_{12}X_{25}X_{34}, &X_{14}X_{22}X_{36}, \\
&X_{12}X_{25}X_{36},& X_{13}X_{25}X_{34}, &X_{13}X_{24}X_{36}, & X_{13}X_{25}X_{36}, & X_{14}X_{25}X_{36} 
 \end{array}
 \end{equation} 
  and Hilbert series (\ref{NoHS36}).  We take $G=[1, 3, 2][5, 4, 6]-[1, 3, 5][2, 4, 6]\in O_F$ and 
  $$G=G_0 \mod(FV)$$ 
  with 
  $$G_0=
  -X_{11}X_{12}X_{23}X_{24}X_{35}X_{36} 
  + X_{11}X_{12}X_{24}X_{25}X_{33}X_{36} 
  + X_{11}X_{14}X_{22}X_{23}X_{35}X_{36}.$$
 Only the second and third monomial in $G_0$ can be used to extend the coherent matching $T$ to a coherent matching of $\cM_3\cup \{G\}$ and these two choices are indeed equivalent up to symmetry. 
The final check is that the coherent matching 
$T,  X_{11}X_{12}X_{24}X_{25}X_{33}X_{36}$  gives the correct Hilbert series.      
 Finally we observe that up to symmetry  the only vertex of Type 2 of $\cM_3$ that does not give a \Sagbi{} basis is the one with only even coordinates. 
 
   {\bf  (Type 3)} 
 They correspond to coherent matchings that do not involve the variables 
 $$\FV=\{X_{21}, X_{32}, X_{33}, X_{14}, X_{15}, X_{25}\} $$
 and hence vertices of the Newton polytope  
  $$ \begin{pmatrix}
  			* & * &*  & 0& 0 &  * \\
  			0 &* &* & *   &0  &  *  \\
  			* &0 & 0 & * &* & * \\
  		\end{pmatrix} 
.$$		
 Here 
\begin{equation}
\label{Type3mat}
 X^{(3)}=\begin{pmatrix}
  			X_{11} &X_{12} &X_{13}   & 0 & 0  &  X_{16}  \\
  			0 &X_{22} &X_{23}  & X_{24} & 0 & X_{26}\\
  			X_{31} &0 & 0 & X_{34} &X_{35} & X_{36} \\
  		\end{pmatrix} 
\end{equation} 	
and the $\New(\cM_3^{(3)})$ has $92$ vertices and $24$ orbits under the subgroup (of order $4$) of $G_{3\times 6}$ fixing the shape above. Only one orbit does not correspond to a \Sagbi{} basis and it has orbit representative:
$$
\begin{pmatrix}
8& 8&2& 0& 0& 2 \\
0& 2&8&8&0&2 \\
2&0&0&2& 10&6
  \end{pmatrix}  
  $$ associated coherent matching 
  \begin{equation} 
\begin{array}{cccccc}
T=&
X_{12}X_{23}X_{31},& X_{12}X_{24}X_{31}, & X_{11}X_{22}X_{35}, & X_{11}X_{22}X_{36},& X_{11}X_{23}X_{34}, \\
&X_{11}X_{23}X_{35}, & X_{11}X_{23}X_{36}, & X_{11}X_{24}X_{35}, & X_{11}X_{24}X_{36},& X_{11}X_{26}X_{35}, \\
&X_{12}X_{23}X_{34},& X_{12}X_{23}X_{35}, & X_{12}X_{23}X_{36}, & X_{12}X_{24}X_{35},& X_{12}X_{24}X_{36}, \\
&X_{12}X_{26}X_{35}, & X_{13}X_{24}X_{35}, & X_{13}X_{24}X_{36}, & X_{16}X_{23}X_{35},& X_{16}X_{24}X_{35}, \end{array}
 \end{equation} 
and Hilbert series (\ref{NoHS36}).  We take $G=[3, 4, 1][2, 5, 6]-[3, 4, 2][1,5,6]\in O_F$ and 
$$
G=G_0 \mod(FV)
$$ 
with 
\begin{multline*}
G_0=
-X_{11}X_{13}X_{22}X_{26}X_{34}X_{35} + 
X_{11}X_{16}X_{22}X_{23}X_{34}X_{35} -\\ 
X_{12}X_{13}X_{24}X_{26}X_{31}X_{35} + 
X_{13}X_{16}X_{22}X_{24}X_{31}X_{35}.
\end{multline*}
Only the second and third monomial in $G_0$ can be used to extend the coherent matching $T$ to a coherent matching of $\cM_3\cup \{G\}$, and these two choices are indeed equivalent up to symmetry. 
The final check is that the coherent matching 
$T, X_{11}X_{16}X_{22}X_{23}X_{34}X_{35}$  gives the correct Hilbert series.      
 Finally we observe that up to symmetry  the only vertex of Type 3 of $\cM_3$ that does not give a \Sagbi{} basis is the one with only even coordinates.

 {\bf  (Type 4)} 
 They correspond to coherent matchings that do not involve the variables 
 $$\FV=\{X_{11}, X_{12}, X_{23}, X_{24}, X_{35}, X_{36}\} $$
 and hence vertices of the Newton polytope  
  $$ \begin{pmatrix}
  			0 & 0 &*  & * & *  &  * \\
  			*  &* & 0 & 0   &*   &  *  \\
  			* & *  & *  & * &0 & 0  \\
  		\end{pmatrix} 
.$$		
 Here 
\begin{equation}
\label{Type4mat}
 X^{(4)}=\begin{pmatrix}
  			0 & 0 & X_{13}   & X_{14}  & X_{15}   &  X_{16}  \\
  			X_{21}  &X_{22} &0  & 0 & X_{25} & X_{26}\\
  			X_{31} &X_{32}  & X_{33}  & X_{34} & 0 & 0 \\
  		\end{pmatrix} 
\end{equation} 	
and the $\New(\cM_3^{(4)})$ has $160$ vertices and $6$ orbits under the subgroup of order $48$ of $G_{3\times 6}$ fixing the shape above. Only one orbit does not correspond to a \Sagbi{} basis and it has orbit representative:
$$
\begin{pmatrix}
0& 0&2 &8& 8& 2  \\
8& 2& 0& 0& 2& 8  \\
2& 8& 8& 2& 0& 0 
  \end{pmatrix}  
  $$ associated coherent matching 
  \begin{equation} 
\begin{array}{cccccc}
T=&X_{13}X_{21}X_{32},& X_{14}X_{21}X_{32},& X_{15}X_{21}X_{32},& X_{16}X_{21}X_{32},& X_{14}X_{21}X_{33},\\
    &X_{15}X_{21}X_{33},& X_{16}X_{21}X_{33},& X_{15}X_{21}X_{34},& X_{14}X_{26}X_{31},& X_{15}X_{26}X_{31},\\  
    &X_{14}X_{22}X_{33},& X_{15}X_{22}X_{33},& X_{13}X_{26}X_{32},& X_{14}X_{25}X_{32},& X_{14}X_{26}X_{32},\\
    &X_{15}X_{26}X_{32},& X_{14}X_{25}X_{33},& X_{14}X_{26}X_{33},& X_{15}X_{26}X_{33},& X_{15}X_{26}X_{34}, 
 \end{array}
 \end{equation} 
  and Hilbert series (\ref{NoHS36}).   We take $G=[2, 3, 4][5, 1, 6]-[2, 3, 5][4, 1, 6] \in O_F$ and 
  $$G=G_0 \mod(FV)$$ 
  with 
 \begin{multline*}
G_0= X_{13}X_{14}X_{25}X_{26}X_{31}X_{32} + X_{13}X_{15}X_{22}X_{26}X_{31}X_{34}\\ + X_{13}X_{16}X_{21}X_{25}X_{32}X_{34} - X_{13}X_{16}X_{22}X_{25}X_{31}X_{34}\\ + X_{14}X_{16}X_{22}X_{25}X_{31}X_{33} + X_{15}X_{16}X_{21}X_{22}X_{33}X_{34}. 
 \end{multline*}
Only the first and last  monomial in $G_0$ can be used to extend the coherent matching $T$ to a coherent matching of $\cM_3\cup \{G\}$ and these two choices are indeed equivalent up to symmetry. 
The final check is that the coherent matching 
$T,  X_{13}X_{14}X_{25}X_{26}X_{31}X_{32}$  gives the correct Hilbert series.      
 Finally we observe that up to symmetry  the only vertex of Type 4 of $\cM_3$ that does not give a \Sagbi{} basis is the one with only even coordinates. 
\end{proof} 
 
 A surprising corollary of the proof of  Theorem \ref{univSagbi3x6} is
 
 \begin{corollary}
 \label{onlyeven}
 The vertices of the Newton polytope of $\cM_3$ that correspond to  \Sagbi{} bases of $G(3,6)$ are exactly those 
 that have at least one odd coordinate. 
 \end{corollary} 
 
Theorem \ref{univSagbi3x6} and its proof provide the basic data for the universal \Sagbi{} basis of $G(3,7)$. Again it is a proof by computation, but it is impossible to list all data in the same detail as in the proof of Theorem \ref{univSagbi3x6}. Let us first extend the family $O_F$. Again we start from
$$
F=[1,2,3][4,5,6]-[1,2,4][3,5,6], 
$$
this time as an element of $G(3,7)$, and let $O_F^{(7)}$ be the orbit, up to sign, of $F$ under the action of $S_7$ by permutation of the columns.  Therefore $O_F^{(7)}$  consists of $7\cdot 15=105$ polynomials  of degree $2$ in $G(3,7)$. 

\begin{theorem}\label{univSagbi3x7}
The union of $\cM_3\subset G(3,7)$ and $O_F^{(7)}$ is a universal \Sagbi{}  basis of $G(3,7)$. More precisely, for every monomial order one obtains a \Sagbi{} basis of $G(3,7)$ by adding to $\cM_3$ at most $3$ elements of $O_F^{(7)}$.
\end{theorem}

Theorem \ref{univSagbi3x7} follows immediately from Proposition \ref{add_elements} that provides precise information on the  nature of elements that must be added to $\cM_3$. Let $G(3,6)^{(i)}$ be the algebra generated by the $3$-minors of the matrix $X_{3\times 7}^{(i)}$ obtained from $X_{3\times 7}$ by removing column $i$ with $i=1,\dots,7$.  

\begin{proposition}\label{add_elements}
Let $A$ be the algebra generated by the monomials of a coherent matching $T$ of  $\cM_3 \subset K[X_{3\times 7}]$. 
 Then the following hold:
\begin{enumerate}
\item Let $h = H(G(3,7),2) - H(A,2)$. Then $h\le 3$.
\item  The set $W_T$ of indices $i$ such that $H(A\cap K[X_{3\times 7}^{(i)}],2)\neq H(G(3,6),2)$   has exactly $h$ elements. 
\item There are  elements $G_i \in  O_F^{(7)}$ with $i\in W_T$ such that every augmentation  of   $T$ to a coherent matching of $\cM_3 \cup \{G_i : i\in W_T\}$ gives an algebra with Hilbert series equal to that of $G(3,7)$. 
\end{enumerate}
\end{proposition}
 
\begin{proof}
(1) This follows from the computation of the Hilbert functions for coherent matchings of $\cM_3$ by SagbiGrass (see Section \ref{tools}). 

(2) We know that $H(G(3,6),2)-H(A\cap K[X_{3\times 7}^{(i)}],2)$ is either $0$ or $1$ and since they correspond to distinct column-degrees  there must be exactly $h$ values of $i$ giving $1$.

(3) For each of the four  types of coherent matching we compute a representative of each orbit. From the exponent matrix $D$  associated to the given matching  $T$ we compute the exponent matrices $D_i$ for the restrictions $T_i$  of $T$ to $K[X_{3\times 7}^{(i)}]$, and check which of the $D_i$ have only even entries. By Corollary \ref{onlyeven} we know that this means that $H(G(3,6),2)-H(A\cap K[X_{3\times 7}^{(i)}],2)=1$.  Hence $W_T$ consists of the $i$ such that $D_i$ has only even entries.  Foreach such $i\in W_T$ the crucial point is to find the suitable element $G_i$. The choice is suggested by a finer analysis of the exponent matrices $D_i$. Looking back at the proof of Theorem \ref{univSagbi3x6} we see that the types of a coherent matching of $G(3,6)$  is given by $4-u_{10}$ where $u_{10}$ counts the entries $10$ in the matrix.  Note that in general the type of $T$ and $T_i$ can differ. Of course, one cannot expect to obtain exactly the representatives of the types as in the proof of Theorem \ref{univSagbi3x6}, and this forces us to substitute the column indices  in the polynomials $G$ described in the proof of Theorem \ref{univSagbi3x6}. After this is done, we have identified the polynomials $G_i \in O_F^{(7)}$ with $i\in W_T$. It remains to check that every augmentation   $T$ to a coherent matching of $\cM_3 \cup \{G_i : i\in W_T\}$ gives an algebra with Hilbert series equal to that of $G(3,7)$ and this is done again with (the extended version of) SagbiGrass. 

Let us give an example. Suppose the given coherent matching $T$ of $\cM_3\subset G(3,7)$ has  exponent matrix  
$$
\begin{pmatrix}
	15 &0 &0&10&2&6 &2\\
	0&15 &0 &4&9&4 &3\\
	0 &0&15 &1&4&5&10\\
\end{pmatrix}
$$ 
of type 1.  We check the restrictions to the various $G(3,6)^{(i)}$ and find out that only for $i=1$ we have a matrix $D_1$ with only even entries, that is $W_T=\{1\}$.  In this case 
$$
D_1=\begin{pmatrix}
0& 0& 0&10&2&6&2\\
0&10& 0& 0&6&2&2\\
0& 0&10& 0&2&2&6
\end{pmatrix} 
$$ which is still of type 1 while the orbit representative for type 1 used in the proof  of Theorem \ref{univSagbi3x6}  is 
 $$
\begin{pmatrix}
10& 0& 0& 6& 2&2 \\ 
0& 10& 0& 2& 6& 2 \\
0& 0&10& 2& 2& 6 
\end{pmatrix} 
$$
 This yields the assignment of columns $1,2,3,4,5,6$ of the type $1$ matrix for $G(3,6)$ to the columns $4,2,3,6,5,7$, respectively,  in the matrix for $G(3,7)$. Then we substitute the column indices in the polynomial
$$
G=[1, 4, 2][5, 3,6]-[1, 4, 5][2, 3,6].
$$
and obtain the polynomial 
$$
G_1= [4,6,2][5,3,7] - [4,6,5][2,3,7].
$$
It remains to find the extensions of the given coherent matching of $\cM_3$ to coherent matchings of $\cM_3\cup \{G_1\}$ and to check their Hilbert functions, which indeed coincide with the Hilbert function of $G(3,7)$.
\end{proof}

Let $<$ be a monomial order on $K[X_{m\times n}]$ with associated coherent match $T$ of $\cM_m$. Let $Z$ be the set of the variables that so not appear in $T$ (the forgotten variables).  We know by \cite{SZ} that $Z$ contains exactly $m(m-1)$ variables, $(m-1)$ per row. The precise combinatorial characterization of the possible $Z$ is given in \cite[Thm.3.6]{SZ}.  
Let $X_{m\times n}^0$ be the matrix obtained form  $X_{m\times n}$ by replacing with $0$ the elements of $Z$.
Let $\phi_Z: K[X_{m\times n}]\to K[X_{m\times n}^0]$ the $K$-algebra map defined by $X_{ij}\to 0$ if $X_{ij}\in Z$ and  $X_{ij}\to X_{ij}$ otherwise. 
 Set  $G(m,n)^0=\phi_Z( G(m,n))$. 
As we have already seen implicitly in the proof of Theorem \ref{univSagbi3x6},  an important consequence of Theorem  \ref{dimok} is the following: 

\begin{corollary} 
\label{forgotten} 
One has  $G(m,n)\simeq G(m,n)^0$ and  $\ini( G(m,n) )=\ini( G(m,n)^0 )$. Furthermore $\cF$ is a \Sagbi{} basis of 
$G(m,n)$ if and only if $\phi_Z(\cF)$ is a \Sagbi{} basis of  $G(m,n)^0$. 
\end{corollary}

\begin{proof} By construction $K[T]\subset \ini(G(m,n)^0)$ and hence $\dim G(m,n)^0\geq \dim K[T]$ and,  by  Theorem \ref{dimok}, 
$\dim G(m,n)=\dim K[T]$. Summing up, we have a surjection of $K$-algebra domains $\phi_Z:G(m,n)\to G(m,n)^0$ and 
$\dim G(m,n)\leq \dim G(m,n)^0$. It follows that  $G(m,n)\simeq G(m,n)^0$.  Now we claim that  the variables in $Z$ cannot be involved in the algebra  $\ini( G(m,n) )$. Suppose, by contradiction,  there is a $X^a$ monomial  in  $\ini( G(m,n))$  that involves some  variables that are in  $Z$. 
Hence $\dim K[T,X^a]=\dim K[T]+1=\dim G(m,n)+1$ and $\dim K[T,X^a]_i\leq \dim \ini( G(m,n))_i=\dim G(m,n)_i$ for every $i\in \NN$. This is a contradiction. For claim it follows that  $\ini(F)=\ini(\phi_Z(F))$ for each $F\in G(m,n)$ and this implies the remaining assertions. \end{proof}

\begin{remark}\label{higher3xN}
Let $\cF$ be a finite set of homogeneous polynomials generating the subalgebra $A$ of the polynomial ring endowed with a monomial order $<$. 
We say that $\cF$ is \emph{non\Sagbi{} in degree $k$} if $K[\ini_<(\cF)]$ and $\ini_<(A)$ coincide in degrees $< k$, but differ in degree $k$. This can of course be controlled by comparing the Hilbert series since $K[\ini_<(\cF)] \subset \ini(A)$.  
 
 As the proof of Theorem \ref{univSagbi3x7} shows, $\cM_3\subset G(3,7)$ has no ``genuine'' coherent non\Sagbi{} matching: if a matching is non\Sagbi, then at least one restriction is non\Sagbi{} for $G(3,6)^{(i)}$. In contrast, $\cM_3\subset G(3,8)$ has many genuine non\Sagbi{} coherent matchings. SagbiGrass computes the complete list of Hilbert series for all coherent matchings of $\cM_3\subset G(3,8)$ in a few hours. The list contains matchings that are non\Sagbi{} in degrees $2,4$ and $5$.  On the other hand, there are no coherent matchings that are non\Sagbi{} in degree $6$ or higher. 
 
While it takes months to compute the full list of Hilbert series for coherent matchings of $\cM_3\subset G(3,9)$ with certainty, one quickly finds coherent matchings that are non\Sagbi{} in degree $6$, and therefore are genuine for $G(3,9)$. Moreover, a random search found the  revlex order given below  such that the   \Sagbi{} basis of $G(3,9)$, computed  by sagbiNormaliz.lib,  is finite but  needs  ``new'' \Sagbi{} elements in degrees $5,8,13$ and $18$ (more precisely,  $12$ elements in degree $5$ and one element in each of the remaining degrees). Indeed we computed a \Sagbi{} basis  and the initial algebra basis of $G(3,9)^0$ and then used Corollary \ref{forgotten}. 
This revlex order is obtained by ordering the variables according to the following table
$$
\begin{pmatrix}
12&6& 18&10&16&24&8& 27&9\\ 
20&2& 26&7& 14&1& 13&11&23\\
15&19&5& 17&4& 25&22&3& 21 
\end{pmatrix}
$$
i.e., the variable corresponding to the number $1$ is the largest, then the variable  corresponding to the number $2$ is the second largest and so on.   
 \end{remark}

 Let $\cF\subset R= K[X_1,\dots,X_n]$ be a set of polynomials of constant degree $m$ and $<$ a monomial order. 
 As a vector space, the algebra $K[\cF]$ generated by $\cF$ is the direct sum of the lowest degree components of the powers $I^k$ of the ideal $I$ generated by $\cF$. The Rees algebra
$$
\Rees(I) = \bigoplus_{k=0}^\infty I^kY^k\subset R[Y]
$$
 has a standard $\ZZ$-grading obtained by extending  the standard grading on $R$ by  $\deg Y = -m+1$ and also a finer standard $\ZZ^2$-grading associated to 
$\deg X_i=(1,0)$ and $\deg Y=(-m,1)$.  The properties of $\Rees(I)$ reflect properties 
of the powers  of $I$.  For example, $I$ has linear powers (i.e., all powers of $I$ have a linear resolution)  if and only if the $(1,0)$-regularity of $\Rees(I)$ is $0$, see \cite{BCV1} or \cite[Chap. 8]{BCRV}. 
Therefore it is interesting to investigate the initial algebras of $\Rees(I)$ under an extension of $<$ to  
$R[Y]$.   Since the elements of $\cF$ have constant degree $m$, $K[\cF]\iso K[\cF Y]$ is an algebra  retract of $\Rees(I)$ by degree selection. As a $K$-algebra $\Rees(I)$ is generated by the indeterminates $X_j$ and $\cF Y$,
and  the critical question is whether they form  a \Sagbi{} basis of $\Rees(I)$.  It is easy to see that necessary 
conditions for this are  that  $\cF$ is a Gr\"obner basis of $I$ and  a \Sagbi{} basis of $K[\cF]$.  Furthermore the variables 
$X_j$ and $\cF Y$ form a \Sagbi{} basis of $\Rees(I)$ if and only if 
 $\cF$ is a Gr\"obner basis of $I$ and 
$$
\ini(\Rees(I)) = \Rees(\ini(I)). 
$$

For $m\leq n$, $\cF= \cM_m \subset K[X_{m\times n}]$ and the ideal $I=(\cM_m)=I_m(X_{m\times n})$ of maximal minors we may   study the various \Sagbi{} bases of $\Rees(m,n)=\Rees(I)$  as we have done for its subalgebra $G(m,n)$. As $\cM_m$ is a universal Gr\"obner basis,  in this case we have that the union of the variables $X_{ij}$ and $\cM_mY$ is a \Sagbi{} basis of $\Rees(m,n)$ with respect to $<$ if and only if $\ini(\Rees(I)) = \Rees(\ini(I))$. 
And this is actually the case for a  diagonal order, see \cite[Thm. 6.2.19]{BCRV}.  

In the investigation of the \Sagbi{} bases of $\Rees(m,n)$ the relevant Newton polytope is simply a translate of $\New(\cM_m)$ since the other algebras generators are the variables. Hence the subdivision  of the coherent matchings in the $4$ types given in  \cite{SZ} and used in the proofs of Theorems \ref{univSagbi3x6} and \ref{univSagbi3x7} applies to the study of  $\Rees(3,n)$ as well.   For a coherent matching $T$ of $\cM_m$ we will say that $T$ is \Sagbi{} for $\Rees(m,n)$ if the variables $X_{ij}$ and $\cM_mY$ form a \Sagbi{} basis for $\Rees(m,n)$ for any monomial order associated to $T$, i.e. $\ini(\Rees(m,n))=\Rees(T)$. Having said this, we do not  try to explore fully this aspect of the story, but  only  collect some data and problems in the following remarks.  

\begin{remark}
\label{R=G}
 Let $T$ be a coherent matching of $\cM_3\subset G(3,6)$.  A computation shows that  $T$ is \Sagbi{} for $\Rees(3,6)$ if and only if $T$ is \Sagbi{} for $G(3,6)$, i.e., (by Corollary \ref{onlyeven}) the  corresponding vertex of $\New(\cM_3)$ has at least one odd entry. In this case it would be interesting to check if the ideal $(T)$ has linear powers. This is in principe a finite check because the $(1,0)$-regularity of $\Rees(T)$ is computable. 
   \end{remark} 

\begin{remark}\label{Rees_1} In view of Remark \ref{R=G} and the proof of   Theorem  \ref{univSagbi3x6}, to understand all the \Sagbi{} bases of $\Rees(3,6)$ it remains to analyze coherent matchings with only even entries. There are exactly four such coherent matchings up to  $G_{3\times 6}$,  one per type. They  are described in the proof of  Theorem  \ref{univSagbi3x6}.  
Let us denote them by $T_1,T_2,T_3,T_4$ where the index denotes the type.    

(a) $T_2,T_3$ and $T_4$ behave in the same way. Their Rees algebras have the same bigraded Hilbert series and are defined by $36$ equations of bidegree $(0,2)$ and $45$ equations of bidegree $(1,1)$. In this cases it suffices to add $GY^2$ (with $G$ the polynomial given in the proof of proof of  Theorem  \ref{univSagbi3x6}) to get a \Sagbi{} basis of $\Rees(3,6)$. 
 
(b) $T_1$ behaves differently. Its Rees algebra has a bigraded Hilbert series that is different from the one of $T_2,T_3,T_4$ and it is defined by $36$ equations of bidegree $(0,2)$, $45$ equations of bidegree $(1,1)$ and  a single equation of bidegree $(3,3)$.  In this case it is not enough to add $GY^2$ to get a \Sagbi{} basis of $\Rees(3,6)$. One needs also to add an element $HY^3$ of bidegree $(3,3)$ to complete the \Sagbi{} basis of $\Rees(3,6)$. The element $H$ is of the form 
 $$-X_{16}X_{24}X_{35}[1,4,5][2,5,6][3,4,6] + X_{15}X_{26}X_{34}[1,4,6][2,4,5][3,5,6].$$

 (c) Note that the data above shows already that the ideal $(T_1)$ does not have linear powers while it is possible that $(T_i)$ has linear powers for $i=2,3,4$. 
 \end{remark} 
  
\begin{remark} \label{Rees_2}
For $\Rees(3,7)$ the nature and structure of the various \Sagbi{} bases  is much more complicated. By SagbiGrass we have computed a complete list of $\ZZ$-graded Hilbert series which shows that type 1 is again the most complicated. Moreover we have done random experiments that use sagbiNormaliz.lib for the computation of \Sagbi{} bases. We content ourselves with two observations that show the high complexity of a universal \Sagbi{} basis for $\Rees(3,7)$.

(a) There are type $1$ coherent matchings $T$ of $\cM_3$ that are \Sagbi{} for $G(3,7)$, but not for $\Rees(3,7)$.

(b) There are coherent matchings $T$ of $\cM_3$  that, to be completed to a \Sagbi{} basis for $\Rees(3,7)$, need several elements in high bidegrees. By a random search similar to the one mentioned at the end of Remark \ref{higher3xN}  we found an example that needs new elements up to bidegree $(12,17)$. At that bidegree, for  lack of RAM,  we stopped the computation.
\end{remark}


\begin{thebibliography}{15}

\bibitem{CoCoA}
J.~Abbott, A.~M. ~Bigatti, L.~Robbiano,  
CoCoA: a system for doing Computations in Commutative Algebra, 
Available at \url{https://cocoa.dima.unige.it} 

\bibitem{AOP}
H.~Abo, G.~Ottaviani, C.~Peterson, 
Non-Defectivity of Grassmannians of planes, 
J. Algebr. Geom., 21 (2012), 1--20.


\bibitem{BZ} 
D.~Bernstein, A.~Zelevinsky, 
Combinatorics of maximal minors, 
J. Algebraic Combin. 2 (1993), no. 2, 111--121. 


\bibitem{Bo} 
A.~Boocher,
{Free resolutions and sparse determinantal ideals},
Math. Res. Lett. 19 (2012), no. 4, 805--821.

\bibitem{BMC}
L.~Bossinger, F.~Mohammadi, A.N.~ Ch\'avez, 
Families of Gr\"obner degenerations, Grassmannians and universal cluster algebras, 
SIGMA 17 (2021). 

\bibitem{BC}
W.~Bruns, A.~Conca, 
KRS and powers of determinantal ideals, 
Compositio Mathematica Volume 111, Issue 1, Pages 111--122, 1998. 

  
\bibitem{BCSagbi}
W.~Bruns, A.~Conca,  
Sagbi combinatorics of maximal minors and a Sagbi algorithm, 
J. Symb. Comput. 120  (2024), Article ID 102237, 14 p.

\bibitem{BCV}
W.~Bruns, A.~Conca, M.~Varbaro, 
Relations between the minors of a generic matrix, 
Adv.~Math.~244 (2013), 171--206.

\bibitem{BCV1}
W.~Bruns, A.~Conca, M.~Varbaro, 
Maximal minors and linear powers, 
J.~Reine Angew.~Math.~702 (2015), 41--53. 
 
\bibitem{BCRV}
W.~Bruns, A.~Conca, C.~Raicu, M.~Varbaro, 
Determinants, Gr\"obner bases and cohomology, 
Springer Monographs in Mathematics, Springer, Cham 2022.

\bibitem{BrGu}
W.~Bruns, J.~Gubeladze, 
Polytopes, rings, and $K$-theory, 
Springer Monographs in Mathematics, Springer, Dordrecht, 2009.

\bibitem{BH}
W.~Bruns,  J.~Herzog, 
Cohen-Macaulay rings, 
Rev.~ed.~Cambridge: Cambridge University Press, 1998.

\bibitem{Nmz}
W.~Bruns, B.~Ichim, C.~S\"oger, U.~von~der~Ohe,	
Normaliz. Algorithms for rational cones and affine monoids, 
Available at \url{https://normaliz.uos.de}.

\bibitem{CS} 
D.~Cartwright, B.~Sturmfels, 
{ The Hilbert scheme of the diagonal in a product of  projective spaces},
Int.~Math.~Res.~Not.~9 (2010), 1741--1771.

\bibitem{CM1}
O.~Clarke, F.~Mohammadi,
Toric degenerations of Grassmannians and Schubert varieties from matching field tableaux, 
Journal of Algebra 559, 646-678 (2020). 

\bibitem{CHM}
O.~Clarke, A.~Higashitani, F.~Mohammadi,
Combinatorial mutations of Gelfand-Tsetlin polytopes, Feigin-Fourier-Littelmann-Vinberg polytopes, and block diagonal matching field polytopes, 
Journal of Pure and Applied Algebra 228 (7), (2024). 

\bibitem{CM2}
O.~Clarke, F.~Mohammadi,
Minimal cellular resolutions of powers of matching field ideals, 
Journal of Pure and Applied Algebra 229 (2025). 


\bibitem{CHV}
A.~Conca, J.~Herzog, G.~Valla,
\Sagbi{} bases with applications to blow-up algebras, 
J.~Reine Angew.~Math.~474 (1996), 113--138.


  \bibitem{C} 
  A.~Conca, 
{ Linear spaces, transversal polymatroids and ASL domains},
J.~Algebraic Combin.~25 (2007), no. 1, 25--41. 


\bibitem{CDG1} 
A.~Conca, E.~De Negri, E.~Gorla,
{ Universal Gr\"obner bases for maximal minors},
Int.~Math.~Res.~Not.~11 (2015), 3245--3262.

\bibitem{CDG2} 
A.~Conca, E.~De Negri, E.~Gorla,
{Universal Gr\"obner bases and Cartwright-Sturmfels ideals},
Int.~Math.~Res.~Not.~7 (2020), 1979--1991.

\bibitem{CDG3} 
A.~Conca, E.~De Negri, E.~Gorla,
{ Multigraded generic initial ideals of determinantal ideals},
Homological and Computational Methods in Commutative Algebra
A.~Conca, J.~Gubeladze, and T.~R\"omer Eds., Springer (2018). 

\bibitem{CDG4} 
A.~Conca, E.~De Negri, E.~Gorla,
{ Cartwright-Sturmfels ideals associated to graphs and linear spaces},
J.~Comb.~Algebra 2, no.~3 (2018), 231--257.

\bibitem{CDG5} 
A.~Conca, E.~De Negri, E.~Gorla,
{Radical generic initial ideals},  
Vietnam J.~Math.~50 (2022). 

\bibitem{CDG6} 
A.~Conca, E.~De Negri, E.~Gorla, 
{Radical support for multigraded ideals,} 
S\`ao Paulo J.~Math. Sci. 17 (2023). 


\bibitem{CW} 
A.~Conca, V.~Welker,
{Lov\'asz-Saks-Schrijver ideals and coordinate sections of determinantal varieties},
 Algebra \& Number Theory 13 (2019), no. 2, 455--484.  



\bibitem{Sing}	
W.~Decker, G.~M.~Greuel, G.~Pfister,  H.~Schonemann, 
Singular 4-1-1 -- A computer algebra system for polynomial computations, 
\url{https://www.singular.uni-kl.de}, 2018.

\bibitem{EneHerz}
V.~Ene, J.~Herzog, 
Gr\"obner bases in commutative algebra, vol. 130 of Graduate Studies
in Mathematics. American Mathematical Society, Providence, RI, 2012.
 

\bibitem{HPPRS}
H.~Huang, M.~Perlman, C.~Polini, C.~Raicu,  A.~Sammartano, 
Relations between the $2\times 2$ minors of a generic matrix, 
 Adv.~Math.~386 (2021). 


\bibitem{KapMad}
D.~Kapur, K.~Madlener, 
A completion procedure for computing a canonical basis for a
$k$-subalgebra,
 In Computers and mathematics (Cambridge, MA, 1989). Springer, New York,
1989, pp. 1--11.

\bibitem{KrRo}
M.~Kreuzer, L.~Robbiano, 
Computational commutative algebra II, 
 Springer-Verlag,  Berlin, 2005.

\bibitem{Lembo}
F.~Lembo,  
SAGBI Bases of Algebras of Minors, 
Master thesis, University of Genova, 2024.  



\bibitem{MS}
F.~Mohammadi, K.~Shaw, 
Toric degenerations of Grassmannians from matching fields, 
Algebraic Combinatorics 2 (6), 1109-1124, (2019). 

\bibitem{M}
F.~Mohammadi, 
Gr\"obner degenerations of determinantal ideals with an application to toric degenerations of Grassmannians, 
International  Congress on Mathematical Software, 285-295, (2024). 


\bibitem{RobSwe}
L.~Robbiano and M.~Sweedler. 
Subalgebra bases. 
In Commutative algebra (Salvador,
1988), vol. 1430 of Lecture Notes in Math.~Springer, Berlin, 1990, pp. 61--87.

\bibitem{Stan}
R.~Stanley, 
Combinatorics and commutative algebra. 
2nd ed. Birkhäuser, Boston 1996.

\bibitem{SZ} 
B.~Sturmfels,  A.~Zelevinsky, 
Maximal minors and their leading terms.
Adv.~Math.~98 (1993), no.~1, 65--112.

\bibitem{StuGrPol}
B.~Sturmfels, 
 Gr\"obner bases and convex polytopes, 
 vol. 8 of University Lecture Series.
American Mathematical Society, Providence, RI, 1996.

\end{thebibliography}
\end{document}